\documentclass[12pt, leqno]{article}
\usepackage{amsfonts}
\usepackage{amssymb,amsmath,latexsym,amscd}
\evensidemargin =\oddsidemargin

\font\teneufm=eufm10
\font\seveneufm=eufm7
\font\fiveeufm=eufm5
\newfam\eufmfam
\textfont\eufmfam=\teneufm
\scriptfont\eufmfam=\seveneufm
\scriptscriptfont\eufmfam=\fiveeufm
\def\frak#1{{\fam\eufmfam\relax#1}}

\def\Char{\text{\rm Char}}

\def\beq{\begin{equation}}
\def\eeq{\end{equation}}
\def\bea{\begin{eqnarray}}
\def\eea{\end{eqnarray}}
\def\beas{\begin{eqnarray*}}
\def\eeas{\end{eqnarray*}}

\def\cplus{\hbox{$\supset${\raise1.05pt\hbox{\kern -0.55em
${\scriptscriptstyle +}$}}\ }}

 \DeclareMathOperator{\Hom}{Hom}

 \DeclareMathOperator{\TL}{TL}
\DeclareMathOperator{\DT}{DT} \DeclareMathOperator{\CTL}{CTL}
\DeclareMathOperator{\CB}{CB} \DeclareMathOperator{\G}{G}

\newtheorem{Def}{Definition}[section]

\newtheorem{Prop}[Def]{Proposition}
\newtheorem{Thm}[Def]{Theorem}
\newtheorem{Lem}[Def]{Lemma}
\newtheorem{Coro}[Def]{Corollary}
\title{Cyclotomic
Temperley-Lieb algebra of type $D$ and its representation theory}

\author{Jie Sun$\hbox{\,}$\vspace{0.3cm}\\$\hbox{\,}$ {\small University of Alberta, Department of Mathematical and Statistical Sciences},\\ {\small 632 CAB, University of Alberta, Edmonton, Alberta, Canada T6G 2G1}\\ {\small Email:\ jsun@math.ualberta.ca}}
\date{}

\begin{document}
\maketitle

\begin{abstract}
We define a new class of algebras, cyclotomic Temperley-Lieb
algebras of type $D$, in a diagrammatic way, which is a
generalization of Temperley-Lieb algebras of type $D$. We prove that
the cyclotomic Temperley-Lieb algebras of type $D$ are cellular. In
fact, an explicit cellular basis is given by means of combinatorial
methods. After determining all the irreducible representations of
these algebras, we give a necessary and sufficient condition for a
cyclotomic Temperley-Lieb algebra of type $D$ to be
quasi-hereditary.
\end{abstract}


\medskip


\bigskip

\section{Introduction}

$\indent$Temperley-Lieb algebras are a class of finite dimensional
 associative algebras first introduced by Temperley and Lieb (1971) in their
analysis of Potts models and later rediscovered by Jones (1983) to
characterize his algebras arising from the tower construction of
semi-simple algebras in the study of subfactors. As well as having
many applications to physics, Temperley-Lieb algebras are also of
great value in several areas of mathematics, including the theory of
quantum groups
 and knot theory, where they are closely related to the Jones polynomial and isotopy
 invariants of links. This relationship was explained in Jones (1987),
 where it was shown that Temperley-Lieb algebras occur
 naturally as quotients of Hecke algebras arising from a
 Coxeter system of type $A$. In his Ph.D. thesis, Graham (1995) studied
 certain quotients of Hecke algebras associated to a Coxeter
 diagram $X$, which were called Temperley-Lieb algebras of type
 $X$. Graham classified finite dimensional Temperley-Lieb
 algebras into seven infinite families: $A$, $B$, $D$, $E$, $F$, $H$ and
 $I$. Some affine versions of Temperley-Lieb algebras have also been studied by Graham and
 Lehrer (1998), Fan and Green (1999), and so on.

Recently, cyclotomic Temperley-Lieb algebras of type $A$, as a
generalization
 of Temperley-Lieb algebras of type $A$, were introduced and investigated by Rui and Xi
 (2004) and were proved to be cellular in the sense of Graham and Lehrer (1996) by means of dotted
 planar graphs. Using the theory of
 cellular algebras (K\"{o}nig and Xi, 1998, 1999a, 1999c) all irreducible
 representations of cyclotomic Temperley-Lieb algebras of type $A$ were parametrised by Rui and Xi (2004). They also determined when these algebras are quasi-hereditary (in the sense of Cline, Parshall and Scott, 1988). Moreover, a necessary and
 sufficient condition for a cyclotomic Temperley-Lieb algebra of type
 $A$ to be semi-simple was provided in Rui, Xi and Yu (2005). Martin and Saleur (1994)
 introduced blob algebras which can be seen as another generalization of the Temperley-Lieb algebras of type
 $A$ by adding an idempotent generator and some
 defining relations. Cao and Zhu (2006) defined cyclotomic
 blob algebras as algebras of diagrams and showed that cyclotomic
 blob algebras are cellular.

 A common approach when studying the cellular structure of the above algebras is via diagrams.
 Diagram calculi have already
 been developed for many algebras: Temperley-Lieb
 algebras of type $A$ (Westbury, 1995, Graham and Lehrer, 1996), Temperley-Lieb
 algebras of type $D$ (Green, 1998), affine Temperley-Lieb
 algebras of type $A$ (Fan and Green, 1999) and cyclotomic Temperley-Lieb
 algebras of type $A$ (Rui and Xi, 2004). The idea behind this is
 Martin and Saleur's pictorial definition of Temperley-Lieb
 algebras of type $A$ (Martin and Saleur, 1994). K\"{o}nig and Wang
 (2008) gave a uniform apporach to cyclotomic extensions of diagram
 algebras.

The purpose of this article is to introduce a cyclotomic version of
Temperley-Lieb algebras of type $D$ by means of diagrammatic
generators and relations. After recalling the definition of
Temperley-Lieb algebras of type $D$ in section 2, we define the
cyclotomic Temperley-Lieb algebras of type $D$ in section 3. We
prove that the cyclotomic Temperley-Lieb algebras of type $D$ are
cellular in section 4 and an explicit cellular basis is given by
combinatorial methods. In section 5, we determine all the
irreducible representations of these algebras and give a necessary
and sufficient condition for a cyclotomic Temperley-Lieb algebra of
type $D$ to be quasi-hereditary. In the last section, we go through
a concrete example.

For simplicity, we suppose that the ground ring $K$ is a field. It
is assumed that all algebras considered in this article are finite
dimensional associative $K$-algebras with identity, all modules are
unitary, and all modules are left modules unless stated otherwise.

 \section{Temperley-Lieb algebras of type $D$}

  $\indent$ In this section we recall the definition of Temperley-Lieb
algebras of type $D$ and introduce the category of $m$-decorated
tangles. The following figure is the Dynkin diagram of type $D_{n}$
which will be used to define Temperley-Lieb algebras of type $D$.

\begin{center}
\setlength{\unitlength}{1cm}
\begin{picture}(10,3)
\put(1,0){\circle*{.25}} \put(1,2){\circle*{.25}}
\put(2,1){\circle*{.25}} \put(3,1){\circle*{.25}}
\put(8,1){\circle*{.25}} \put(9,1){\circle*{.25}}
\put(1,2){\line(1,-1){1}} \put(1,0){\line(1,1){1}}
\put(2,1){\line(1,0){1}} \put(8,1){\line(1,0){1}}
\put(3,1){\dashbox{0.2}(5,0)[]{}} \put(0.5,0){$\overline{1}$}
\put(0.5,2){1} \put(2,0.5){2} \put(3,0.5){3} \put(7.8,0.5){n-2}
\put(9,0.5){n-1}
\end{picture}
\end{center}

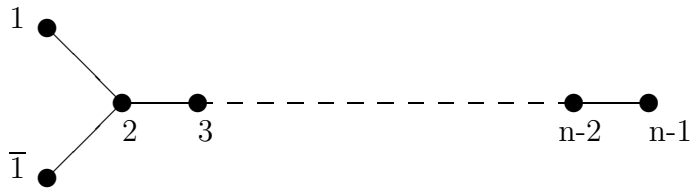
\begin{center}Figure 1. Dynkin diagram of type $D_{n}$
\end{center}

\begin{Def} {\rm Let $n \in \mathbb{N} \geq 4$ and $\delta \in K$ be a
parameter. The \emph{Temperley-Lieb algebra of type} $D$, denoted by
$\TL_{n}(D)$, is an associative algebra over $K$ with generators
$E_{\overline{1}}$, $E_{1}$, $E_{2}$, $\cdots$, $E_{n-1}$ subject to
the following relations:

$\mbox{~~~~~~~~~~~}E_{i}^{2}=\delta E_{i} \mbox { for all } i, $

$\mbox{~~~~~~~~~}E_{i}E_{j}=E_{j}E_{i} \mbox{ if i and j are not
connected in the Dynkin diagram},$

$\mbox{~~~~~~}E_{i}E_{j}E_{i}=E_{i} \mbox{ otherwise}.$}

\label{D2.1.1}
\end{Def}

The approach to define an algebra by diagrams has been used to
understand otherwise purely abstract algebraic objects such as
representations and cellular structures (in the sense of Graham and
Lehrer, 1996). The category of decorated tangles was introduced in
Green (1998) to study Temperley-Lieb algebras of types $B$ and $D$.
The idea behind this is Martin and Saleur's pictorial definition of
Temperley-Lieb algebras of type $A$ (see Martin and Saleur, 1994).
We will recall the basic notions of tangles and decorated tangles in
Green (1998) for later use.

\begin{Def} {\rm Let $p$ and $q$ be positive integers. A
\emph{tangle} of type $(p,q)$ is a portion of a knot diagram
contained in a rectangle in the plane, consisting of arcs and closed
cycles, such that the endpoints of the arcs consist of $p$ points in
the top edge of the rectangle and $q$ points in the bottom edge.}
\end{Def}

For simplicity, a tangle of type $(n,n)$ is said to be a tangle of
type $n$. We refer to the boundary of a rectangle as its
\emph{frame}. Two tangles are equal if there exists an isotopy of
the plane carrying one to the other such that the two diagrams can
be identified when we fix the frame of each rectangle (see Freyd and
Yetter, 1989 or Kauffman, 1990 for details). The endpoints of arcs
are called \emph{vertices}. The vertices in the top (respectively,
bottom) edge of the frame are numbered consecutively starting with
vertex number $1$ at the leftmost end. Arcs in a diagram are called
\emph{horizontal arcs} if they connect two vertices sitting in the
same edge of the frame; and are called \emph{vertical arcs} if they
connect two vertices sitting in the different edges.

\begin{Def} {\rm A \emph{decorated tangle} is a
crossing-free tangle in which each arc is assigned a nonnegative
integer. Any arc or closed cycle not exposed to the left face of the
frame (namely, this arc or closed cycle is separated from the left
face of the frame by another arc) must be assigned the integer
zero.}
\end{Def}

If an arc (respectively, a closed cycle) is assigned the value $s$,
we represent this pictorially by decorating the arc (respectively,
the closed cycle) with $s$ blobs. Each blob is marked by a hollow
disc. If there are too many blobs sitting in an arc (or a closed
cycle), we will mark the arc (or the closed cycle) with a
representative blob, and write the value $s$ around the blob. For
any horizontal arc that links vertices $i$ and $j$ with $i<j$, we
denote it by $\{ i<j \}$. For any vertical arc that links vertices
$i$ in the top edge and $j$ in the bottom edge, we denote it by $\{
i, j \}$.

Before defining the category of decorated tangles, we first need to
give the following rules for movements of blobs between connected
arcs:

(1) A blob in an arc can move freely to another arc if the two arcs
share a common endpoint.

(2) In the above movement, the number of blobs is given by adding
the number of blobs from each arc.

\begin{Def} {\rm The \emph{category $\DT$ of decorated tangles}
is defined as follows:

(1) The objects of $\DT$ are positive integers.

(2) The morphisms from $p$ to $q$ are the decorated tangles of type
$(p,q)$.

(3) For any $G_{1} \in \Hom_{\DT}(p,q)$ and $G_{2} \in
\Hom_{\DT}(q,r)$, the composition $G_{1} \diamond G_{2}$ is defined
to be the concatenation of the tangle $G_{1}$ above the tangle
$G_{2}$, identifying the bottom vertices of $G_{1}$ with the top
vertices of $G_{2}$ and assigning the nonnegative integers of arcs
and closed cycles according to the rules for movements of blobs
between connected arcs.}
\end{Def}

{\it Remark}. Note that for there to be any morphisms from $p$ to
$q$, it is necessary for $p+q$ to be even. A careful calculation in
Green (1998) shows that $(G_{1} \diamond G_{2})\diamond G_{3}= G_{1}
\diamond (G_{2}\diamond G_{3})$. The category-theoretic definition
allows us to define an algebra of decorated tangles.

\begin{Def} {\rm Let $n \in \mathbb{N}$. The algebra $\DT_{n}$
has a $K$-basis consisting of morphisms from $n$ to $n$ in $\DT$ and
the multiplication is given by the composition in $\DT$.}
\end{Def}

It is convenient to define certain named decorated tangles,
$e_{\bar{1}}$, $e_{1}$, $e_{2}$, $\cdots$, $e_{n-2}$ and $e_{n-1}$
in the algebra $\DT_{n}$. Define $e_{i}$, $1 \leq i \leq n-1$, to be
the decorated tangle in which both the top edge and the bottom edge
have a horizontal arc $\{ i < i+1 \}$, and the other arcs are
vertical. There are no blobs on any of the arcs in $e_{i}$ (for $1
\leq i \leq n-1$). Define $e_{\bar{1}}$ to be the decorated tangle
obtained from $e_{1}$ by adding one blob to each of the two
horizontal arcs. Now we can realize the Temperley-Lieb algebra of
type $D$ in terms of decorated tangles, and this is due to Green
(1998).

\begin{Thm}{\rm (Green, 1998)}. Let $n \in \mathbb{N} \geq 4$ and $\delta \in K$ be a
parameter. The algebra $\widetilde{\TL}_{n}(D)$ has a $K$-basis $\{
e_{\bar{1}}, e_{1}, e_{2}, \cdots, e_{n-1} \}$ and the
multiplication is induced from that of ${\DT}_{n}$ subject to the
following relations:
\begin{center}
\setlength{\unitlength}{1cm}
\begin{picture}(8,8)
\put(3,7){\oval(1.5,0.8)} \put(4,7){=} \put(5,7){$\delta$}
\put(3,3.5){\line(0,1){2.5}} \put(3,4.5){\circle{.25}}
\put(3,5.3){\circle{.25}} \put(4,4.9){=}
\put(5,3.5){\line(0,1){2.5}} \put(1.5,1.5){\oval(1.5,0.8)}
\put(3.5,0.3){\line(0,1){2.5}} \put(4,1.5){=}
\put(5.5,1.5){\oval(1.5,0.8)} \put(7,0.3){\line(0,1){2.5}}
\put(2.2,1.5){\circle{.25}} \put(3.5,1.5){\circle{.25}}
\put(6.3,1.5){\circle{.25}}
\end{picture}
\end{center}

\begin{center}{\rm Figure 2. Relations for $\widetilde{\TL}_{n}(D)$}
\end{center}

There is a basis for $\widetilde{\TL}_{n}(D)$ which is in natural
bijection with elements of $\DT_{n}$ which have at most one blob on
each arc or closed cycle, and which satisfy one of the following two
mutually exclusive conditions:

(I) The diagram contains one closed cycle on which there is a blob,
and there are no other closed cycles with blobs in the diagram.
Also, there is at least one horizontal arc in the diagram.

(II) The diagram contains no closed cycles and the total number of
blobs is even.

We say that an element of $\DT_{n}$ which satisfies these hypotheses
is \emph{$D$-admissible diagram }of type I or type II, depending on
which of the two conditions above it satisfies.

There is an isomorphism $\rho_{D}: \TL_{n}(D) \longrightarrow
\widetilde{\TL}_{n}(D)$ which takes $E_{\overline{1}}$ to
$e_{\bar{1}}$ and $E_{i}$ to $e_{i}$ for all $1 \leq i \leq n-1$.
\end{Thm}

The first relation implies that any closed cycle with no blobs can
be removed and the resulting diagram multiplied by parameter
$\delta$. The second relation gives that any arc or closed cycle
with $r$ (for $r>1$) blobs is equivalent to the arc or closed cycle
which carries $r-2$ blobs. The third relation yields that any arc or
closed cycle loses its blob in the presence of a closed cycle with
one blob. Using the first and third relations, all closed cycles may
be removed from the resulting diagram except the last closed cycle
with a blob. The following lemma gives the dimension of
$\TL_{n}(D)$.

\begin{Lem}{\rm (Green, 1998)}. Let $C(n)$ be the Catalan number $C(n):= \frac{1}{n+1} \binom{2n}{n}$. In $\widetilde{\TL}_{n} (D)$, the number of $D$-admissible
diagrams of type I is $C(n)-1$, and the number of $D$-admissible
diagrams of type II is $\frac{1}{2} \binom{2n}{n}$. This is a total
of $(\frac{n+3}{2}) C(n)-1$ which is the dimension of $\TL_{n}(D)$.
\end{Lem}

To define the cyclotomic version of Temperley-Lieb algebras of type
$D$, we need the notion of $m$-decorated tangles.

\begin{Def} {\rm Let $p$, $q$, and $m$ be positive integers. An
$m$\emph{-decorated tangle} of type $(p,q)$ is a crossing-free
tangle of type $(p,q)$ in which each arc (and each closed cycle, if
any) is assigned a pair of nonnegative integers $[r,s]$ such that
$r$ is at most $m-1$, and such that $s$ is assigned zero if the arc
(or the closed cycle) is not exposed to the left face of the frame.}
\end{Def}

{\it Remark}. If an arc (respectively, a closed cycle) is assigned
the value $[r,s]$, we represent this pictorially by decorating the
arc (respectively, the closed cycle) with $r$ dots and $s$ blobs.
Each dot is marked by a filled disc, while each blob is marked by a
hollow disc. If there are too many dots or blobs sitting on an arc
(or a closed cycle), we will mark the arc (or the closed cycle) with
a representative dot or a representative blob, and write the value
$[r,s]$ around the dot or the blob. In the case $m=1$, a
$m$-decorated tangle is a decorated tangle as defined before.

From now on, we make the following convention: Let $\{i<j\}$ be a
horizontal arc, and assume that there are some dots or blobs sitting
on the arc. There is no need to distinguish where the blobs sit, but
we must distinguish the dots sitting at the left side from those at
the right side of the arc. We call $i$ the \emph{left endpoint} and
call $j$ the \emph{right endpoint} of the given arc $\{i<j\}$, and
the dots sitting at the left (respectively, right) endpoint are
called \emph{left} (respectively, \emph{right}) \emph{dots} of the
arc. We always assume that in an $m$-decorated tangle all dots are
left dots. In terms of vertical arcs and closed cycles we do not
define their left endpoints and right endpoints.

To define the category of $m$-decorated tangles, we need first to
give the following rules for movements of dots and blobs between
connected arcs:

(1) A blob and a dot can interchange.

(2) A left dot of a horizontal arc $\{i<j\}$ is equal to $m-1$ right
dots  of the arc $\{i<j\}$, and conversely, a right dot of the
horizontal arc is equal to $m-1$ left dots.

(3) A blob in an arc can move freely to another arc if the two arcs
share a common endpoint. In the movement, the number of blobs is
given by adding the number of blobs from each arc.

(4) A dot in a vertical arc $\{i,j\}$ can move to another vertical
arc if the two arcs share a common endpoint. In the movement, the
number of dots are given directly by sum.

(5) A right dot of a horizontal arc $\{i<j\}$ can move to another
horizontal arc $\{j<k\}$ (or $\{h<j\}$), and this dot will be
considered as a left dot of the arc $\{j<k\}$ (or a right dot of the
arc $\{h<j\}$). Similarly, a left dot of a horizontal arc $\{i<j\}$
can move to another horizontal arc $\{k<i\}$ (or $\{i<h\}$), and
this dot will be considered as a right dot of the arc $\{k<i\}$ (or
a left dot of the arc $\{i<h\}$).

(6) A dot in a vertical arc $\{i,j\}$ can move to a horizontal arc
$\{k<i\}$ or $\{k<j\}$, and this dot will be regarded as a right dot
on the arc $\{k<i\}$ or $\{k<j\}$. A dot in a vertical arc $\{i,j\}$
can also move to a horizontal arc $\{k>i\}$ or $\{k>j\}$, and this
dot will be regarded as a left dot on the arc $\{k>i\}$ or
$\{k>j\}$.

(7) A left dot of a horizontal arc $\{i<j\}$ can move to a vertical
arc $\{i,k\}$ or $\{k,i\}$. Similarly, A right dot of a horizontal
arc $\{i<j\}$ can move to a vertical arc $\{j,k\}$ or $\{k,j\}$.

In the above movements, numbers of dots are reduced modulo $m$.

\begin{Def} {\rm Let $m \in \mathbb{N}$. The \emph{category $\DT_{m}$ of $m$-decorated
tangles }is defined as follows:

(1) The objects of $\DT_{m}$ are positive integers.

(2) The morphisms from $p$ to $q$ are the $m$-decorated tangles of
type $(p,q)$.

(3) For any $G_{1}\in\Hom_{\DT_{m}}(p,q)$ and
$G_{2}\in\Hom_{\DT_{m}}(q,r)$, the composition $G_{1}\diamond G_{2}$
is defined to be the concatenation of the tangle $G_{1}$ above the
tangle $G_{2}$, identifying the bottom vertices of $G_{1}$ with the
top vertices of $G_{2}$ and assigning the nonnegative integer pairs
$[r,s]$ of arcs and closed cycles according to the rules for
movements of dots and blobs between connected arcs.}\label{D2.2.2}
\end{Def}

By a careful calculation we can check that $(G_{1} \diamond
G_{2})\diamond G_{3}= G_{1} \diamond (G_{2}\diamond G_{3})$
according to the rules for movements of dots and blobs defined as
before. In the case $m=1$, the category $\DT_{m}$ of $1$-decorated
tangles is exactly the category of decorated tangles. We end this
subsection by defining an algebra of $m$-decorated tangles based on
the previous category-theoretic definition.

\begin{Def} {\rm Let $m,n \in \mathbb{N}$. The algebra $\DT_{m,n}$
has a $K$-basis consisting of morphisms from $n$ to $n$ in
$\DT_{m}$, and the multiplication is given by the composition in
$\DT_{m}$.}
\end{Def}

\section{Cyclotomic
Temperley-Lieb algebras of type $D$}

$\indent$ In this section we define a new class of algebras called
cyclotomic Temperley-Lieb algebras of type $D$ via diagrams. This is
motivated by the work of Rui and Xi (2004), Cao and Zhu (2006) and
Green (1998). We shall first focus on some special $m$-decorated
tangles, which play a key role in describing the diagram calculus
relevant to cyclotomic Temperley-Lieb algebras of type $D$.

\begin{Def} {\rm Let $n \in \mathbb{N} \geq 4$ and $m\in\mathbb{N}$. An \emph{$m$-cyclotomic $D$-admissible
diagram} of type $n$ is a $m$-decorated tangle which has at most one
blob on each arc or closed cycle, and which satisfy one of the
following two mutually exclusive conditions:

(I) The diagram contains one closed cycle on which there is a blob
and no dots, and no other closed cycles or blobs in the diagram.
Also there is at least one horizontal arc in the diagram.

(II) The diagram contains no closed cycles and the total number of
blobs is even.

We say that an $m$-cyclotomic $D$-admissible diagram is of type I or
type II, depending on which of the two conditions above it
satisfies.}\label{D3.1.1}
\end{Def}

The following figure shows a typical example of $2$-cyclotomic
$D$-admissible diagrams of type $n=4$.
\begin{center}
\setlength{\unitlength}{1cm}
\begin{picture}(10,3)
\put(2.5,0){type I} \put(8,0){type II} \put(1,1.75){\oval(1,0.75)}
\put(1,2.1){\circle{.2}} \put(4.5,1.5){\circle*{.2}}
\put(2,1){\line(1,1){2}} \put(3,1){\line(1,1){2}}
\put(9,1){\line(0,1){2}} \put(10,1){\line(0,1){2}}
\put(2.5,3){\oval(1,1)[b]} \put(4.5,1){\oval(1,1)[t]}
\put(7.5,3){\oval(1,1)[b]} \put(7.5,1){\oval(1,1)[t]}
\put(2.5,2.5){\circle*{.2}} \put(7.2,2.6){\circle*{.2}}
\put(7.8,2.6){\circle{.2}} \put(9,2){\circle{.2}}
\put(7.5,1.5){\circle*{.2}} \put(10,2){\circle*{.2}}
\end{picture}
\end{center}

\begin{center}Figure 3. $2$-Cyclotomic
$D$-admissible diagrams of type $n=4$
\end{center}

With the definition of $m$-cyclotomic $D$-admissible diagram, we
define cyclotomic Temperley-Lieb algebras of type $D$ as follows.

\begin{Def}
{\rm Let $n \in \mathbb{N} \geq 4$ and $m\in\mathbb{N}$. Let $K$ be
a field, and let $\delta_{i} \in K$ ($0\leq i\leq m-1$) with
$\delta_{i}\delta_{0}=\delta_{i}$ for $1\leq i\leq m-1$. The
\emph{cyclotomic Temperley-Lieb algebra of type $D$}, denoted by
$\CTL(D)_{m,n}$, has a $K$-basis consisting of $m$-cyclotomic
$D$-admissible diagrams of type $n$ and the multiplication is
induced from that of ${\DT}_{m,n}$ subject to the following
relations:}
\begin{center}
\setlength{\unitlength}{1cm}
\begin{picture}(8,8)
\put(3,7){\oval(1.5,0.8)} \put(3,7.4){\circle*{.25}} \put(4,7){=}
\put(5,7){$\delta_{i}$} \put(2.7,7.8){[i,0]}
\put(3,3.5){\line(0,1){2.5}} \put(3,4.5){\circle{.25}}
\put(3,5.3){\circle*{.25}} \put(5,4.8){\circle*{.25}} \put(4,4.9){=}
\put(5,3.5){\line(0,1){2.5}} \put(1.8,4.8){[i,2]}
\put(5.5,4.8){[i,0]}

\put(1.5,1.5){\oval(1.5,0.8)} \put(3.5,0.3){\line(0,1){2.5}}
\put(4,1.5){=} \put(5.5,1.5){\oval(1.5,0.8)}
\put(7,0.3){\line(0,1){2.5}} \put(2.2,1.5){\circle{.25}}
\put(3.5,1.2){\circle{.25}} \put(6.3,1.5){\circle{.25}}
\put(7,1.5){\circle*{.25}} \put(3.5,1.9){\circle*{.25}}
\put(1.2,2.3){[0,1]} \put(5.2,2.3){[0,1]} \put(3.6,2.3){[i,1]}
\put(7.5,2.3){[i,0]}
\end{picture}
\end{center}
\begin{center}
\setlength{\unitlength}{1cm}
\begin{picture}(10,1)
\put(2.5,1){\oval(1.5,0.8)} \put(1.75,1){\circle*{.25}}
\put(3.25,1){\circle{.25}} \put(4,1){=} \put(5,1){$\delta_{i}$}
\put(2.2,1.8){[i,1]} \put(6.5,1){\oval(1.5,0.8)}
\put(7.25,1){\circle{.25}} \put(6.2,1.8){[0,1]} \put(8,1){($i \geq
1$)}
\end{picture}
\end{center}\label{D3.2.1}
\end{Def}
\begin{center}Figure 4. Relations for $\CTL(D)_{m,n}$
\end{center}

The first relation gives that any closed cycle with $i$ dots and no
blobs can be removed with the resulting diagram being multiplied by
the parameter $\delta_{i}$ to compensate. The second relation
implies that any arc or closed cycle with $i$ dots and $r$ (for
$r>1$) blobs is equivalent to the arc or closed cycle which carries
$i$ dots and $r-2$ blobs. The third relation yields that any arc or
closed cycle loses its blob in the presence of a closed cycle with
one blob. The fourth relation implies that any closed cycle with $i$
(for $i \geq 1$) dots and one blob loses its dots with the resulting
diagram being multiplied by the parameter $\delta_{i}$ to
compensate. Thus, if we denote by $G_{1}\circ G_{2}$ the diagram
induced from $G_{1}\diamond G_{2}$ (see Definition \ref{D2.2.2})
according to the relations in Definition \ref{D3.2.1}, we can give
explicitly the expression of the multiplication:
$$G_{1}\cdot
G_{2}=(\prod_{i=1}^{m-\!1}\delta_{i}^{n(\bar{i},G_{1},G_{2})^{+}}\prod_{i=0}^{m-\!1}\delta_{i}^{n(\bar{i},G_{1},G_{2})^{-}})G_{1}\circ
G_{2},$$ where $n(\bar{i},G_{1},G_{2})^{+}$ is the number of closed
cycles with $i$ dots and one blob in $G_{1}\diamond G_{2}$ (see
Definition \ref{D2.2.2}), and $n(\bar{i},G_{1},G_{2})^{-}$ is the
number of closed cycles with $i$ dots and no blobs in $G_{1}\diamond
G_{2}$. We can check that $(G_{1}\cdot G_{2})\cdot G_{3}=G_{1}\cdot
(G_{2}\cdot G_{3})$. In fact $$(G_{1}\cdot G_{2})\cdot
G_{3}=(\delta_{0}^{N_{0}}\prod_{i=1}^{m-\!1}\delta_{i}^{N_{i}})(G_{1}\circ
G_{2})\circ G_{3},$$
$$G_{1}\cdot (G_{2}\cdot
G_{3})=(\delta_{0}^{N_{0}^{'}}\prod_{i=1}^{m-\!1}\delta_{i}^{N_{i}^{'}})G_{1}\circ
(G_{2}\circ G_{3}),$$ where
$N_{0}=n(\bar{0},G_{1},G_{2})^{-}+n(\bar{0},G_{1}\circ
G_{2},G_{3})^{-}$,
$N_{i}=n(\bar{i},G_{1},G_{2})^{+}+n(\bar{i},G_{1},G_{2})^{-}+n(\bar{i},G_{1}\circ
G_{2},G_{3})^{+}+n(\bar{i},G_{1}\circ G_{2},G_{3})^{-}$ and
$N_{0}^{'}=n(\bar{0},G_{2},G_{3})^{-}+n(\bar{0},G_{1}, G_{2}\circ
G_{3})^{-}$,
$N_{i}^{'}=n(\bar{i},G_{2},G_{3})^{+}+n(\bar{i},G_{2},G_{3})^{-}+n(\bar{i},G_{1},
G_{2}\circ G_{3})^{+}+n(\bar{i},G_{1}, G_{2}\circ G_{3})^{-}$. Since
$(G_{1}\circ G_{2})\circ G_{3}$ is induced from $(G_{1} \diamond
G_{2})\diamond G_{3}$ and $G_{1}\circ (G_{2}\circ G_{3})$ is induced
from $G_{1} \diamond (G_{2}\diamond G_{3})$ according to the
relations in Definition \ref{D3.2.1}, we know that $(G_{1}\circ
G_{2})\circ G_{3}=G_{1}\circ (G_{2}\circ G_{3})$ from $(G_{1}
\diamond G_{2})\diamond G_{3}= G_{1} \diamond (G_{2}\diamond
G_{3})$. Note that the total number of closed cycles with $i$ (for
$i\geq 1$) dots and one blob and closed cycles with $i$ (for $i\geq
1$) dots and no blobs in $G_{1}\diamond G_{2}$ and $(G_{1}\circ
G_{2})\diamond G_{3}$ is equal to the total number of closed cycles
with $i$ (for $i\geq 1$) dots and one blob and closed cycles with
$i$ (for $i\geq 1$) dots and no blobs in $G_{2}\diamond G_{3}$ and
$G_{1}\diamond (G_{2}\circ G_{3})$. Thus $N_{i}=N_{i}^{'}$ (for
$1\leq i\leq m-1$). If there exists $N_{i}\neq 0$ for some $1\leq
i\leq m-1$, then
$\delta_{0}^{N_{0}}\prod_{i=1}^{m-\!1}\delta_{i}^{N_{i}}=\delta_{0}^{N_{0}^{'}}\prod_{i=1}^{m-\!1}\delta_{i}^{N_{i}^{'}}$
from $\delta_{i}\delta_{0}=\delta_{i}$ for $1\leq i\leq m-1$. If
not, because the total number of closed cycles with no dots and one
blob and closed cycles with no dots and no blobs in $G_{1}\diamond
G_{2}$ and $(G_{1}\circ G_{2})\diamond G_{3}$ is equal to the total
number of closed cycles with no dots and one blob and closed cycles
with no dots and no blobs in $G_{2}\diamond G_{3}$ and
$G_{1}\diamond (G_{2}\circ G_{3})$, we have
$\delta_{0}^{N_{0}}=\delta_{0}^{N_{0}^{'}}$.

Cyclotomic Temperlely-Lieb algebras of type $D$ are a generalization
of the usual Temperley-Lieb algebras of type $D$ (see Graham, 1995,
Green, 1998) on the one hand and have cyclotomic Temperley-Lieb
algebras of type $A$ (see Rui and Xi, 2004) with parameter
$\delta_{0}=1$ as a class of subalgebras on the other hand. If
$m=1$, then $\CTL(D)_{m,n}$ is the usual Temperley-Lieb algebra of
type $D$ with $K$-dimension $(\frac{n+3}{2}) C(n)-1$. From this, it
is clear that the cyclotomic Temperley-Lieb algebra of type $D$ has
$K$-dimension $m^{n}((\frac{n+3}{2}) C(n)-1)$. Note that cyclotomic
Temperlely-Lieb algebras of type $D$ are very different from
cyclotomic blob algebras $\CB_{m,n}$ (see Cao and Zhu, 2006 for the
definition) not only in their basis diagrams but also in their
parameters and relations for multiplication. The diagrams in
$\CTL(D)_{m,n}$ may have a closed circle with one blob while the
diagrams in $\CB_{m,n}$ must not contain any closed circles. And the
total number of blobs in $m$-cyclotomic $D$-admissible diagrams of
type II must be even, while the diagrams in $\CB_{m,n}$ do not have
this restriction. In addition, $\CTL(D)_{m,n}$ have $m$ parameters
while $\CB_{m,n}$ have $2m$ parameters.

It is convenient for us to introduce the following notations for
later use. Denote by $Q(n,k)$ the set of all $m$-cyclotomic
$D$-admissible diagrams in which both the top edge and the bottom
edge have $n$ vertices and $k$ horizontal arcs, by $Q^{+}(n,k)$
$(1\leq k \leq [\frac{n}{2}])$ the subset of $Q(n,k)$ that consists
of $m$-cyclotomic $D$-admissible diarams of type I and by
$Q^{-}(n,k)$ $(0 \leq k \leq [\frac{n}{2}])$ the subset of $Q(n,k)$
that consists of $m$-cyclotomic $D$-admissible diagrams of type II.
In the case $k=\frac{n}{2}$, we define $Q_{1}^{-}(n,\frac{n}{2})$
(respectively, $Q_{2}^{-}(n,\frac{n}{2})$) to be the subset of
$Q^{-}(n,\frac{n}{2})$ that consists of diagrams which have even
(respectively, odd) blobs both in the upper part and in the lower
part of the diagrams. We define $P^{+}(n,k)$ (respectively,
$P^{-}(n,k)$) to be the vector space spanned by all $m$-cyclotomic
$D$-admissible diagrams in $Q^{+}(n,k)$ (respectively,
$Q^{-}(n,k)$). In the case $k=\frac{n}{2}$, we define
$P_{1}^{-}(n,\frac{n}{2})$ (respectively,
$P_{2}^{-}(n,\frac{n}{2})$) to be the vector space spanned by all
$m$-cyclotomic $D$-admissible diagrams in $Q_{1}^{-}(n,\frac{n}{2})$
(respectively, $Q_{2}^{-}(n,\frac{n}{2})$). We also define certain
named $m$-cyclotomic $D$-admissible diagrams
$\mathcal{T}_{1},\cdots,\mathcal{T}_{n}$ in the algebra
$\CTL(D)_{m,n}$.

Define $\mathcal{T}_{i}$ (for $1\leq i\leq n$) to be the
$m$-cyclotomic $D$-admissible diagram in which all arcs are vertical
with no blobs, the $i$-th vertical arc carries one dot, and there
are no other dots.

In the case $n=4$, the $m$-cyclotomic $D$-admissible diagrams
$\mathcal{T}_{1}$, $\mathcal{T}_{2}$, $\mathcal{T}_{3}$ and
$\mathcal{T}_{4}$ are as shown in the following figure.
\begin{center}
\setlength{\unitlength}{1cm}
\begin{picture}(8,3)
\put(1.5,0){$\mathcal{T}_{1}$} \put(6.5,0){$\mathcal{T}_{2}$}
\put(0,1){\line(0,1){2}} \put(1,1){\line(0,1){2}}
\put(2,1){\line(0,1){2}} \put(3,1){\line(0,1){2}}
\put(5,1){\line(0,1){2}} \put(6,1){\line(0,1){2}}
\put(7,1){\line(0,1){2}} \put(8,1){\line(0,1){2}}
\put(0,2){\circle*{.2}} \put(6,2){\circle*{.2}}
\end{picture}
\end{center}

\begin{center} \setlength{\unitlength}{1cm}
\begin{picture}(8,3)
\put(1.5,0){$\mathcal{T}_{3}$} \put(6.5,0){$\mathcal{T}_{4}$}
\put(0,1){\line(0,1){2}} \put(1,1){\line(0,1){2}}
\put(2,1){\line(0,1){2}} \put(3,1){\line(0,1){2}}
\put(5,1){\line(0,1){2}} \put(6,1){\line(0,1){2}}
\put(7,1){\line(0,1){2}} \put(8,1){\line(0,1){2}}
\put(2,2){\circle*{.2}} \put(8,2){\circle*{.2}}
\end{picture}
\end{center}
\begin{center}Figure 5. $m$-Cyclotomic $D$-admissible diagrams
$\mathcal{T}_{1}$, $\mathcal{T}_{2}$, $\mathcal{T}_{3}$ and
$\mathcal{T}_{4}$
\end{center}

\section{Cellular structure of
 $\CTL(D)_{m,n}$}

 $\indent$In this section, we
 investigate the cellular structure of cyclotomic Temperley-Lieb
 algebras of type $D$. After recalling the definition of cellular algebras, we
 construct a cellular basis for $\CTL(D)_{m,n}$ using
 combinatorial methods.

\begin{Def} {\rm (Graham and Lehrer, 1996). An associative $K$-algebra $A$ is called
a \emph{cellular algebra} with cell datum $(\Lambda,M,C,i )$ if the
following conditions are satisfied:

$(C1)$ The finite set $\Lambda$ is partially ordered. Associated
with each $\lambda\in\Lambda$ there is a finite indexing set
$M(\lambda)$. The algebra $A$ has a $K$-basis $C_{S,T}^{\lambda}$
where $(S,T)$ runs through all elements of $M(\lambda)\times
M(\lambda)$ for all $\lambda\in\Lambda$.

 $(C2)$ The map
$i$ is a $K$-linear anti-automorphism of $A$ with $i^{2}=id$ which
sends $C_{S,T}^{\lambda}$ to $C_{T,S}^{\lambda}$ for all
$\lambda\in\Lambda$ and all $S$ and $T$ in $M(\lambda)$.

$(C3)$ For each $\lambda\in\Lambda$ and $S,\ T\in M(\lambda)$ and
each $a\in A$, the product $aC_{S,T}^{\lambda}$ can be written as
$(\sum_{U\in{M(\lambda)}}r_{a}(U,S)C_{U,T}^{\lambda})+r'$, where
$r'$ is a linear combination of basis elements with upper index
$\mu$ strictly smaller than $\lambda$, and where the coefficients
$r_{a}(U,S)\in K$ are independent of $T$.}\label{D4.1.1}
\end{Def}

The basis $\{C_{S,\,T}^{\lambda}|\,\lambda\in\Lambda$ and $S,\,T\in
M(\lambda)\}$ satisfying the above condition is called a
\emph{cellular basis} of $A$. The $K$-linear anti-automorphism $i$
of $A$ with $i^{2}=id$ is called an \emph{involution}. Whether an
algebra is cellular or not depends on the choice of the involution
$i$. A cellular algebra can have more than one cellular basis and
both the poset $\Lambda$ and the indexing sets $M(\lambda)$ may vary
dramatically between different cellular bases of the same algebra.
The size of the poset $\Lambda$ can also be different for different
cellular bases of an algebra (see K\"{o}nig and Xi, 1999b). We now
recall the ring theoretic definition of cellular algebras, which is
equivalent to the original one.

\begin{Def}{\rm (K\"onig and Xi, 1998). Let $A$ be a $K$-algebra
with an involution $i$. A two-sided ideal $J$ of $A$ is called a
\emph{cell ideal} if and only if $i(J)=J$ and there exists a left
ideal $W\subset J$ such that there is an isomorphism of
$A$-bimodules $\alpha : J\simeq W\!\!\otimes_{K}\!i(W)$ making the
following diagram commutative:
$$\CD J@>\alpha>>\!W\!\!\otimes_{K}\!i(W)\\
@ViVV@VVx\otimes y \mapsto i(y)\otimes i(x)V\\
J@>\alpha>>\!W\!\!\otimes_{K}\!i(W)
\endCD
$$
The algebra $A$ $($with the involution $i$$)$ is called
\emph{cellular} if and only if there is a vector space decomposition
$A=J'_{1}\oplus J'_{2}\oplus \cdots \oplus J'_{n}$ $($for some
$n$$)$ with $i(J'_{j})=J'_{j}$ for each $j$ and such that setting
$J_{j}=\oplus_{l=1}^{j}J'_{l}$ gives a chain of two-sided ideals of
$A: 0=J_{0}\subset J_{1}\subset J_{2}\subset \cdots \subset J_{n}=A$
$($each of them fixed by $i$$)$ and for each $j$ $(j=1,\cdots,n)$
the quotient $J'_{j}=J_{j}/J_{j-\!1}$ is a cell ideal $($with
respect to the involution induced by $i$ on the quotient$)$ of
$A/J_{j-\!1}$. We call this chain a \emph{cell chain} for the
cellular algebra $A$.}\label{D4.1.2}
\end{Def}

{\it Remark}. It has been shown that the class of cellular algebras
includes a number of well-known algebras such as Ariki-Koike
algebras (see Ariki and Koike, 1994), Brauer algebras (see Graham
and Lehrer, 1996, K\"{o}nig and Xi, 1999a), Jones' annular algebras
(see Graham and Lehrer, 1996), and partition algebras (see Xi, 1999)
as well as Birman-Wenzl algebras (see Xi, 2000).

\bigskip

We recall a simple example of cellular algebras in Rui and Xi (2004)
for later use. Take $\G_{m,n}$ to be the $K$-subalgebra of the
cyclotomic Temperley-Lieb algebra of type $D$ generated by
$\mathcal{T}_{1},\cdots,\mathcal{T}_{n}$. $\G_{m,n}$ is isomorphic
to the group algebra of the abelian group
$\oplus_{i=1}^{n}\mathbb{Z}/m\mathbb{Z}$. Assume that $K$ is a
splitting field of $x^{m}-1$, then we can write
$\mathcal{T}_{i}^{m}=1$ as
$\prod_{j=1}^{m}(\mathcal{T}_{i}-\xi_{j})=0$ for some
$\xi_{1},\cdots,\xi_{m}\in K$. Let
$\Lambda(m,n)=\{(i_{1},\cdots,i_{n})\,|\,1\leq i_{k}\leq m\}$ for
$n\geq1$, and assume that in the case $n=0$ the set $\Lambda(m,n)$
consists of only one element $\emptyset$. We define
$(i_{1},\cdots,i_{n})\preceq(j_{1},\cdots, j_{n})$ if and only if
$i_{k}\leq j_{k}$ for all $1\leq k\leq n$. For each
$I=(i_{1},\cdots,i_{n})$, define
$C_{1,1}^{I}=\prod_{j=1}^{n}\prod_{l=i_{j}+\!1}^{m}(\mathcal{T}_{j}-\xi_{l})$,
where the product over empty set is assumed to be $1$. We observe
that $\{C_{1,1}^{I}\,|\,I\in\Lambda(m,n)\}$ is a cellular basis of
the algebra $\G_{m,n}$ with respect to the identity involution. In
order to construct a cellular basis for $\CTL(D)_{m,n}$ we introduce
two further notions.

\begin{Def}
{\rm An \emph{$m$-decorated dangle} of type $(n,k)$ is a
crossing-free diagram consisting of $n$ vertices $\{1,\cdots,n\}$,
$k$ horizontal arcs, $n-2k$ vertical lines and some closed cycles
which satisfy the following conditions.

(1) Each horizontal arc (and each closed cycle, if any) carries at
most $m-1$ dots, and each vertical line does not carry any dots.

(2) Only the leftmost vertical line and the horizontal arcs (and the
closed cycles, if any) which appear to the left of the leftmost
vertical line and to the outermost of any nested horizontal arcs may
carry at most one blob.}
\end{Def}

\begin{Def} {\rm An \emph{$m$-cyclotomic $D$-admissible
dangle} of type $(n,k)$ is an $m$-decorated dangle of type $(n,k)$
which satisfies one of the following two mutually exclusive
conditions:

(I) The diagram contains one closed cycle on which there is a blob
and there are no dots, and no other closed cycles or blobs in the
diagram. Also there is at least one horizontal arc in the diagram.

(II) The diagram contains no closed cycles and the total number of
blobs is even if $k\neq\frac{n}{2}$.

We say that an $m$-cyclotomic $D$-admissible dangle is of type I or
type II, depending on which of the two conditions above it
satisfies.}
\end{Def}

We denote by $D(n,k)$ the set of all $m$-cyclotomic $D$-admissible
dangles of type $(n,k)$, by $D^{+}(n,k)$ $(1\leq k \leq
[\frac{n}{2}])$ the subset of $D(n,k)$ that consists of
$m$-cyclotomic $D$-admissible dangles of type I and by $D^{-}(n,k)$
$(0 \leq k \leq [\frac{n}{2}])$ the subset of $D(n,k)$ that consists
of $m$-cyclotomic $D$-admissible dangles of type II. In the case
$k=\frac{n}{2}$, we define $D_{1}^{-}(n,\frac{n}{2})$ (respectively,
$D_{2}^{-}(n,\frac{n}{2})$) to be the subset of
$D^{-}(n,\frac{n}{2})$ that consists of dangles with an even
(respectively, odd) number of blobs in total. We define $V^{+}(n,k)$
(respectively, $V^{-}(n,k)$) to be the vector space spanned by all
$m$-cyclotomic $D$-admissible dangles in $D^{+}(n,k)$ (respectively,
$D^{-}(n,k)$). In the case $k=\frac{n}{2}$, we define
$V_{1}^{-}(n,\frac{n}{2})$ (respectively,
$V_{2}^{-}(n,\frac{n}{2})$) to be the vector space spanned by all
$m$-cyclotomic $D$-admissible dangles in $D_{1}^{-}(n,\frac{n}{2})$
(respectively, $D_{2}^{-}(n,\frac{n}{2})$). The following lemma
shows a close relationship between $P^{+}(n,k)$ (respectively,
$P^{-}(n,k)$) and $V^{+}(n,k)$ (respectively, $V^{-}(n,k)$).

\begin{Lem}
There are three $K$-module isomorphisms:

(1) $P^{+}(n,k)\simeq
V^{+}(n,k)\otimes_{K}V^{+}(n,k)\otimes_{K}G_{m,n-\!2k}$ for $1\leq k
\leq [\frac{n}{2}]$.

(2) $P^{-}(n,k)\simeq
V^{-}(n,k)\otimes_{K}V^{-}(n,k)\otimes_{K}G_{m,n-\!2k}$ for $1\leq k
\leq [\frac{n}{2}]$ and $k\neq\frac{n}{2}$.

(3) $P^{-}_{i}(n,k)\simeq
V^{-}_{i}(n,k)\otimes_{K}V^{-}_{i}(n,k)\otimes_{K}G_{m,0}$
$(i=1,2)$. \label{L4.2.3}
\end{Lem}
{\it Proof.} We first prove (1). Suppose that $G\in P^{+}(n,k)$.
Then $G$ has $n-2k$ vertical arcs. Let $m_{i}$ be the number of dots
in the $i$-th vertical arc and let
$\omega=\mathcal{T}_{1}^{m_{1}}\cdots
\mathcal{T}_{n-\!2k}^{m_{n-\!2k}}$. Clearly, we have $\omega\in
G_{m,n-\!2k}$. Cutting all vertical arcs and omitting all dots in
the vertical arcs, we can divide the diagram $G$ into two half
diagrams, each with a closed cycle carrying one blob. Denote by
$G_{1}$ the upper part, and by $G_{2}$ the lower part. Consequently,
we can write $G$ uniquely as $G_{1}\otimes G_{2}\otimes\omega$,
which belongs to
$V^{+}(n,k)\otimes_{K}V^{+}(n,k)\otimes_{K}G_{m,n-\!2k}$.
Conversely, given such an expression $G_{1}\otimes
G_{2}\otimes\omega$, we have a unique diagram $G$ with $G_{1}$ on
the top and $G_{2}$ on the bottom, the ends being joined in the
unique way which creates a crossing-free diagram. Omitting one
closed cycle in $G$, we get a $m$-cyclotomic $D$-admissible diagram
of type I in $P^{+}(n,k)$. This proves (1).

The proof of (2) is similar to that of (1). The difference is that
before cutting the leftmost vertical arc we first omit the blob (if
any) on it. If the total number of blobs in the horizontal arcs of
$G_{1}$ (respectively, $G_{2}$) is odd, then the leftmost vertical
line will be assigned one blob such that $G_{1}$ (respectively,
$G_{2}$) belongs to $V^{-}(n,k)$.

In case of $k=\frac{n}{2}$, there are no vertical arcs in the
$m$-cyclotomic $D$-admissible diagrams. Thus the proof of (3) is
straightforward.\ \hfill $\square$

\bigskip

As a result, we have the following equivalent basis of
$\CTL(D)_{m,n}$. However, this basis is usually not a cellular
basis. Denote by $\overline{G}_{m,n-\!2k}$ the set of a $K$-basis of
$G_{m,n-\!2k}$. Let $V=\{v_{1}\otimes
v_{2}\otimes\omega\,|\,v_{1},v_{2}\in D^{+}(n,k),\ \omega\in
\overline{G}_{m,n-\!2k},\ 1\leq
k\leq[\frac{n}{2}]\}\cup\{v_{1}\otimes
v_{2}\otimes\omega\,|\,v_{1},v_{2}\in D^{-}(n,k),\ \omega\in
\overline{G}_{m,n-\!2k},\ 0\leq k\leq[\frac{n}{2}]\ and \
k\neq\frac{n}{2}\}$ and  $V'=\{v_{1}\otimes
v_{2}\otimes1\,|\,v_{1},v_{2}\in
D^{-}_{1}(n,\frac{n}{2})\}\cup\{v_{1}\otimes
v_{2}\otimes1\,|\,v_{1},v_{2}\in D^{-}_{2}(n,\frac{n}{2})\}$.

\begin{Coro}The set $V$ (respectively, $V\cup
V'$) constitutes a basis for $\CTL(D)_{m,n}$ if $n$ is odd
(respectively, even). \label{C4.2.4}
\end{Coro}

Now we describe the cellular structure of $\CTL(D)_{m,n}$. Recall
that $\Lambda(m,n-\!2k)=\{(i_{1},\cdots,i_{n-\!2k})\,|\,1\leq
i_{j}\leq m\}$. If $n$ is odd, let
$\Lambda_{m,n}=\{(\widetilde{k},I)^{+},\ (k,J)^{-}\,|\,1\leq
\widetilde{k}\leq[\frac{n}{2}],\ 0\leq k\leq[\frac{n}{2}],\ I\in
\Lambda(m,n-\!2\widetilde{k}),\ J\in \Lambda(m,n-\!2k)\}$. We define
a partial order ``\,$\leqq$\,'' on $\Lambda_{m,n}$:

(1) $(\widetilde{k},I)^{+}\ \leqq \ (k,J)^{-}$ for any
$\widetilde{k}, k, I, J$;

(2) $(\widetilde{k}_{1},I_{1})^{+}\ \leqq \
(\widetilde{k}_{2},I_{2})^{+}$ if $\widetilde{k}_{1}>
\widetilde{k}_{2}$;

(3) $(\widetilde{k},I_{1})^{+}\ \leqq \ (\widetilde{k},I_{2})^{+}$
if $I_{1}\preceq I_{2}$;

(4) $(k_{1},J_{1})^{-}\ \leqq \ (k_{2},J_{2})^{-}$ if $k_{1}>
k_{2}$;

(5) $(k,J_{1})^{-}\ \leqq \ (k,J_{2})^{-}$ if $J_{1}\preceq
J_{2}$;\\From the relations (2) and (3) we know that
$\{(\widetilde{k},I)^{+}\,|\,1\leq \widetilde{k}\leq[\frac{n}{2}],\
I\in \Lambda(m,n-\!2\widetilde{k})\}$ is a partially ordered set, as
is $\{(k,J)^{-}\,|\,0\leq k\leq[\frac{n}{2}],\ J\in
\Lambda(m,n-\!2k)\}$ from relations (4) and (5). Note that
$\Lambda_{m,n}$ is in fact a disjoint union of these two subsets
with the former strictly smaller than the latter by relation (1).
From the partial orderings of $\{(\widetilde{k},I)^{+}\,|\,1\leq
\widetilde{k}\leq[\frac{n}{2}],\ I\in
\Lambda(m,n-\!2\widetilde{k})\}$ and $\{(k,J)^{-}\,|\,0\leq
k\leq[\frac{n}{2}],\ J\in \Lambda(m,n-\!2k)\}$ we can check that
$\leqq$ is reflexive, transitive and antisymmetric, thus
$(\Lambda_{m,n},\leqq)$ is a finite partially ordered set. If $n$ is
even, let $\Lambda_{m,n}=\{(\widetilde{k},I)^{+},\
(\frac{n}{2},\emptyset)_{1}^{-},\ (\frac{n}{2},\emptyset)_{2}^{-} ,\
(k,J)^{-}\,|\,1\leq \widetilde{k}\leq[\frac{n}{2}],\ 0\leq
k\leq[\frac{n}{2}]-1,\ I\in \Lambda(m,n-\!2\widetilde{k}),\ J\in
\Lambda(m,n-\!2k)\}$. In this case, we substitute
$(\widetilde{k},I)^{+}\ \leqq \ (\frac{n}{2},\emptyset)_{i}^{-} \
(i=1,2)\ \leqq\ (k,J)^{-}$ for relation (1) above and retain the
four other relations to define a partial order on $\Lambda_{m,n}$.
Similarly, we can prove that $(\Lambda_{m,n},\leqq)$ is a finite
partially ordered set.

For each $(\widetilde{k},I)^{+}\in\Lambda_{m,n}$, define
$M((\widetilde{k},I)^{+})=\{(v,I)\,|\,v\in
D^{+}(n,\widetilde{k})\}$. Similarly, define
$M((k,J)^{-})=\{(v,J)\,|\,v\in D^{-}(n,k)\}$ for each
$(k,J)^{-}\in\Lambda_{m,n}$. In the case $k=\frac{n}{2}$, define
$M((\frac{n}{2},\emptyset)_{i}^{-})=\{(v,\emptyset)\,|\,v\in
D_{i}^{-}(n,k)\}\ (i=1,2)$.

In the following we use the cellular bases
$\{C^{I}_{1,1}\,|\,I\in\Lambda(m,n-\!2\widetilde{k})\}$ of
$G_{m,n-\!2\widetilde{k}}$ and
$\{C^{J}_{1,1}\,|\,J\in\Lambda(m,n-\!2k)\}$ of $G_{m,n-\!2k}$ to
construct a cellular basis for $\CTL(D)_{m,n}$. For each
$(\widetilde{k},I)^{+}\in\Lambda_{m,n}$ and $v_{1}, v_{2}\in
D^{+}(n,\widetilde{k})$, let
$C^{(\widetilde{k},I)^{+}}_{(v_{1},I),(v_{2},I)}=v_{1}\otimes
v_{2}\otimes C^{I}_{1,1}$. Similarly, let
$C^{(k,J)^{-}}_{(v_{1},J),(v_{2},J)}=v_{1}\otimes v_{2}\otimes
C^{J}_{1,1}$ for each $(k,J)^{-}\in\Lambda_{m,n}$ and $v_{1},
v_{2}\in D^{-}(n,k)$. In the case $k=\frac{n}{2}$, let
$C^{(\frac{n}{2},\emptyset)^{-}_{i}}_{(v_{1},\emptyset),(v_{2},\emptyset)}=v_{1}\otimes
v_{2}\otimes C^{\emptyset}_{1,1}$ for each $v_{1}, v_{2}\in
D^{-}_{i}(n,\frac{n}{2})\ (i=1,2)$. Let
$$C(\widetilde{k},I)^{+}=\left\{C^{(\widetilde{k},I)^{+}}_{(v_{1},I),(v_{2},I)}\,|\,v_{1},
v_{2}\in D^{+}(n,\widetilde{k})\right\}$$
$$C(k,J)^{-}=\left\{C^{(k,J)^{-}}_{(v_{1},J),(v_{2},J)}\,|\,v_{1},
v_{2}\in D^{-}(n,k)\right\}$$
$$C\left(\frac{n}{2},\emptyset\right)^{-}_{i}=\left\{C^{(\frac{n}{2},\emptyset)^{-}_{i}}_{(v_{1},\emptyset),(v_{2},\emptyset)}\,|\,v_{1},
v_{2}\in D^{-}_{i}(n,\frac{n}{2})\right\}\ (i=1,2)$$\\
By Corollary \ref{C4.2.4}, the set $$\biggl(\
\bigcup_{\widetilde{k}=1}^{[\frac{n}{2}]}\ \bigcup_{I\in
\Lambda(m,n-\!2\widetilde{k})}C(\widetilde{k},I)^{+}\ \biggr)\
\bigcup\ \biggl(\ \bigcup_{k=0}^{[\frac{n}{2}]}\ \bigcup_{J\in
\Lambda(m,n-\!2k)}C(k,J)^{-}\ \biggr)$$\\ forms a basis for
$\CTL(D)_{m,n}$ if $n$ is odd and the set
$$\biggl(\bigcup_{\widetilde{k}=1}^{[\frac{n}{2}]}\ \bigcup_{I\in
\Lambda(m,n-\!2\widetilde{k})}C(\widetilde{k},I)^{+}\
\biggr)\bigcup\ \biggl(\bigcup_{k=0}^{[\frac{n}{2}]-1}\
\bigcup_{J\in
\Lambda(m,n-\!2k)}C(k,J)^{-}\ \biggr)\bigcup\ \biggl(\bigcup_{i=1}^{2}C(\frac{n}{2},\emptyset)_{i}^{-}\biggr)$$\\
forms a basis for $\CTL(D)_{m,n}$ if $n$ is even.

The involution $i$ corresponds to top-bottom inversion of an
$m$-cyclotomic $D$-admissible diagram. The following lemma shows in
detail how $i$ acts on $\CTL(D)_{m,n}$.
\begin{Lem}
Let $i$ be the $K$-linear anti-automorphism of $\CTL(D)_{m,n}$ which
corresponds to top-bottom inversion of an $m$-cyclotomic
$D$-admissible diagram. Then we have the following.

(1) $i$ sends $v_{1}\otimes v_{2}\otimes C^{I}_{1,1}$ to
$v_{2}\otimes v_{1}\otimes C^{I}_{1,1}$ for all $v_{1}, v_{2}\in
D^{+}(n,\widetilde{k})$ and $C^{I}_{1,1}\in
\overline{G}_{m,n-\!2\widetilde{k}}$.

(2) $i$ sends $v_{1}\otimes v_{2}\otimes C^{J}_{1,1}$ to
$v_{2}\otimes v_{1}\otimes C^{J}_{1,1}$ for all $v_{1}, v_{2}\in
D^{-}(n,k)$ and $C^{J}_{1,1}\in \overline{G}_{m,n-\!2k}$.

(3) $i$ sends $v_{1}\otimes v_{2}\otimes C^{\emptyset}_{1,1}$ to
$v_{2}\otimes v_{1}\otimes C^{\emptyset}_{1,1}$ for all $v_{1},
v_{2}\in D^{-}_{i}(n,\frac{n}{2})\ (i=1,2)$ if $k=\frac{n}{2}$.
\end{Lem}
{\it Proof.} We first prove (1). Suppose $G=v_{1}\otimes
v_{2}\otimes C^{I}_{1,1}\in C(\widetilde{k},I)^{+}$ which has a
closed cycle with one blob, $\widetilde{k}$ horizontal arcs
$\{i_{11}<i_{12} \}$, $\{i_{21}<i_{22} \}$, $\cdots$,
$\{i_{\widetilde{k}1}<i_{\widetilde{k}2} \}$ in the top edge,
$\widetilde{k}$ horizontal arcs $\{j_{11}<j_{12} \}$,
$\{j_{21}<j_{22} \}$, $\cdots$,
$\{j_{\widetilde{k}1}<j_{\widetilde{k}2} \}$ in the bottom edge and
$n-\!2\widetilde{k}$ vertical arcs $\{i_{1},j_{1} \}$,
$\{i_{2},j_{2} \}$, $\cdots$, $\{i_{\widetilde{k}},j_{\widetilde{k}}
\}$. Then $v_{1}$ is the $m$-cyclotomic $D$-admissible dangle of
type I with $\widetilde{k}$ horizontal arcs $\{i_{11}<i_{12} \}$,
$\{i_{21}<i_{22} \}$, $\cdots$,
$\{i_{\widetilde{k}1}<i_{\widetilde{k}2} \}$ and
$n-\!2\widetilde{k}$ vertical lines $i_{1}$, $i_{2}$, $\cdots$,
$i_{\widetilde{k}}$, and $v_{2}$ is the $m$-cyclotomic
$D$-admissible dangle of type I with $\widetilde{k}$ horizontal arcs
$\{j_{11}<j_{12} \}$, $\{j_{21}<j_{22} \}$, $\cdots$,
$\{j_{\widetilde{k}1}<j_{\widetilde{k}2} \}$ and
$n-\!2\widetilde{k}$ vertical lines $j_{1}$, $j_{2}$, $\cdots$,
$j_{\widetilde{k}}$. Let $G^{'}$ be the top-bottom inversion of $G$.
We can describe $G^{'}$ explicitly, that is, $G^{'}$ has a closed
cycle with one blob, $\widetilde{k}$ horizontal arcs
$\{j_{11}<j_{12} \}$, $\{j_{21}<j_{22} \}$, $\cdots$,
$\{j_{\widetilde{k}1}<j_{\widetilde{k}2} \}$ in the top edge,
$\widetilde{k}$ horizontal arcs $\{i_{11}<i_{12} \}$,
$\{i_{21}<i_{22} \}$, $\cdots$,
$\{i_{\widetilde{k}1}<i_{\widetilde{k}2} \}$ in the bottom edge and
$n-\!2\widetilde{k}$ vertical arcs $\{j_{1},i_{1} \}$,
$\{j_{2},i_{2} \}$, $\cdots$, $\{j_{\widetilde{k}},i_{\widetilde{k}}
\}$. Using the cutting method in the proof of Lemma \ref{L4.2.3} we
know that $G^{'}=v_{2}\otimes v_{1}\otimes C^{I}_{1,1}$.

The proof of (2) is more complicated. Suppose $G=v_{1}\otimes
v_{2}\otimes C^{J}_{1,1}\in C(k,J)^{-}$, where $v_{1}$ is the
$m$-cyclotomic $D$-admissible dangle of type II with $k$ horizontal
arcs $\{i_{11}<i_{12} \}$, $\{i_{21}<i_{22} \}$, $\cdots$,
$\{i_{k1}<i_{k2} \}$ and $n-\!2k$ vertical lines $i_{1}$, $i_{2}$,
$\cdots$, $i_{k}$ and $v_{2}$ is the $m$-cyclotomic $D$-admissible
dangle of type II with $k$ horizontal arcs $\{j_{11}<j_{12} \}$,
$\{j_{21}<j_{22} \}$, $\cdots$, $\{j_{k1}<j_{k2} \}$ and $n-\!2k$
vertical lines $j_{1}$, $j_{2}$, $\cdots$, $j_{k}$. Suppose the
total number of blobs in $v_{1}$ is $2s$ and the total number of
blobs in $v_{2}$ is $2t$, where $s,t\in \mathbb{N}\cup\{ 0 \}$. Then
there are two cases for $v_{1}$: One is that there is a blob on the
vertical line $i_{1}$ and the total number of  blobs on the
horizontal arcs is $2s-1$ ($s\geq 1$), and the other is that there
are no blobs on the vertical line $i_{1}$ and the total number of
blobs on the horizontal arcs is $2s$. Similarly there are two cases
for $v_{2}$. So there are four cases of $G$ for us to consider. In
the case $v_{1}$ has a blob on the vertical line $i_{1}$ and $2s-1$
blobs on its horizontal arcs and $v_{2}$ has a blob on the vertical
line $j_{1}$ and $2t-1$ blobs on its horizontal arcs, there are no
blobs on the leftmost vertical arc of $G$ and the total number of
blobs on the horizontal arcs is $2(s+t-1)$. Let $G^{'}$ be the
top-bottom inversion of $G$. Then $G^{'}$ has no blobs on the
leftmost vertical arc $\{j_{1},i_{1}\}$, $2s-1$ blobs on the
horizontal arcs $\{j_{11}<j_{12} \}$, $\{j_{21}<j_{22} \}$,
$\cdots$, $\{j_{k1}<j_{k2} \}$ in the top edge and $2t-1$ blobs on
the horizontal arcs $\{i_{11}<i_{12} \}$, $\{i_{21}<i_{22} \}$,
$\cdots$, $\{i_{k1}<i_{k2} \}$ in the bottom edge. Using the cutting
method in the proof of Lemma \ref{L4.2.3} we know that
$G^{'}=v_{2}\otimes v_{1}\otimes C^{I}_{1,1}$. The other three cases
can be shown similarly. So $i$ sends $v_{1}\otimes v_{2}\otimes
C^{J}_{1,1}$ to $v_{2}\otimes v_{1}\otimes C^{J}_{1,1}$ for all
$v_{1}, v_{2}\in D^{-}(n,k)$ and $C^{J}_{1,1}\in
\overline{G}_{m,n-\!2k}$.

In case of $k=\frac{n}{2}$, the proof is straightforward since there
is no vertical arc in the $m$-cyclotomic $D$-admissible diagrams.\
\hfill $\square$

\bigskip

The following theorem shows that the datum $(\Lambda_{m,n},M,C,i)$
makes the algebra $\CTL(D)_{m,n}$ into a cellular algebra.

\begin{Thm}
Let $K$ be a splitting field of $x^{m}-1$. Then the cyclotomic
Temperley-Lieb algebra of type $D$ over $K$ is a cellular algebra
with cell datum $(\Lambda_{m,n},M,C,i)$.
\end{Thm}
{\it Proof.} By the above construction of the datum
$(\Lambda_{m,n},M,C,i)$, it is clear that the first two conditions
of Definition \ref{D4.1.1} are satisfied. Now, we will verify the
condition $(C3)$ of the definition. We consider the product $a\cdot
C_{S,T}^{\lambda}$ for each $\lambda\in \Lambda_{m,n}$ and $S, T\in
M(\lambda)$ and each $a\in \CTL(D)_{m,n}$. Note that all
$m$-cyclotomic $D$-admissible diagrams of type $n$ form a free
$K$-basis of $\CTL(D)_{m,n}$, so we only have to consider the case
when $a$ is a basis element. We denote by
$\CTL(D)_{m,n}^{\leqq\lambda}$ (respectively,
$\CTL(D)_{m,n}^{<\lambda}$, $\CTL(D)_{m,n}^{\lambda}$) the
$K$-subspace of $\CTL(D)_{m,n}$ spanning by all basis elements with
upper index smaller than (respectively, strictly smaller than, equal
to) $\lambda$.

When $n$ is odd, we have the following four cases to prove.

(1) For each $a=B=B_{1}\otimes B_{2}\otimes\omega\in Q^{+}(n,k)$ and
each $C_{S,T}^{\lambda}=G=G_{1}\otimes G_{2}\otimes C^{I}_{1,1}\in
C(\widetilde{k},I)^{+}$, let $y=B\cdot G$. We need to consider the
horizontal arcs in $B$ and $G$. It is immediate that $y\in
\CTL(D)_{m,n}^{<(\widetilde{k},I)^{+}}$ when $k>\widetilde{k}$ since
there are at least $k$ horizontal arcs in each edge in $y$. The case
$k\leq\widetilde{k}$ is more subtle. It is clear that $y\in
\CTL(D)_{m,n}^{\leqq(\widetilde{k},I)^{+}}$. Suppose $$y\equiv
(\prod_{i=1}^{m-\!1}\delta_{i}^{n(\bar{i},B,G)^{+}}\prod_{i=0}^{m-\!1}\delta_{i}^{n(\bar{i},B,G)^{-}})
B^{'}\otimes G_{2}\otimes\omega^{'}C_{1,1}^{I}$$ (mod
$\CTL(D)_{m,n}^{<(\widetilde{k},I)^{+}}$) $\in
\CTL(D)_{m,n}^{(\widetilde{k},I)^{+}}$, where $B^{'}\in
D^{+}(n,\widetilde{k}),\ \omega^{'}\in G_{m,n-2\!\widetilde{k}}$. We
observe that the eliminated closed cycles are completely determined
by the horizontal arcs in $B_{2}$ and $G_{1}$. Hence the
coefficients $\delta_{i}^{n(\bar{i},B,G)^{+}}$ and
$\delta_{i}^{n(\bar{i},B,G)^{-}}$ depend only on $B_{2}$ and $G_{1}$
and the closed cycles in $B$ and $G$. Note that $B^{'}$ is
determined only by $B_{1}$, $B_{2}$ and $G_{1}$, and $\omega^{'}$
depends on $B_{2}$, $\omega$, $G_{1}$ and $C_{1,1}^{I}$, thus
$B^{'}$ and $\omega^{'}$ are independent of $G_{2}$. Write
$\omega^{'}=\prod_{j=1}^{n-\!2\widetilde{k}}\mathcal{T}_{j}^{k_{j}}$
for some $0\leq k_{j}\leq m-1$, $1\leq j\leq n-2\widetilde{k}$. By a
careful calculation, we know that $\omega^{'}
C_{1,1}^{I}\equiv\prod_{j=1}^{n-\!2\widetilde{k}}\xi_{i_{j}}^{k_{j}}C_{1,1}^{I}$
(mod $G_{m,n-\!2k}^{\prec I}$), where $G_{m,n-\!2k}^{\prec I}$ is
the $K$-subspace of $G_{m,n-2\!\widetilde{k}}$ spanned by
$C_{1,1}^{J}$ with $J$ strictly smaller than $I$. Note also that the
coefficient $\prod_{j=1}^{n-\!2\widetilde{k}}\xi_{i_{j}}^{k_{j}}$ is
independent of $G_{2}$. Therefore,
$$y\equiv(\prod_{i=1}^{m-\!1}\delta_{i}^{n(\bar{i},B,G)^{+}}\prod_{i=0}^{m-\!1}\delta_{i}^{n(\bar{i},B,G)^{-}}
\prod_{j=1}^{n-\!2\widetilde{k}}\xi_{i_{j}}^{k_{j}})B^{'}\otimes
G_{2}\otimes C_{1,1}^{I}$$ (mod
$\CTL(D)_{m,n}^{<(\widetilde{k},I)^{+}}$). By the above argument,
both $B^{'}$ and the coefficient of $B^{'}\otimes G_{2}\otimes
C_{1,1}^{I}$ are independent of $G_{2}$.

(2) For each $a=B=B_{1}\otimes B_{2}\otimes\omega\in Q^{-}(n,k)$ and
each $C_{S,T}^{\lambda}=G=G_{1}\otimes G_{2}\otimes C^{I}_{1,1}\in
C(\widetilde{k},I)^{+}$, let $y=B\cdot G$. It is clear that $y\in
\CTL(D)_{m,n}^{<(\widetilde{k},I)^{+}}$ when $k>\widetilde{k}$. In
the case $k\leq \widetilde{k}$, we have $y\in
\CTL(D)_{m,n}^{\leqq(\widetilde{k},I)^{+}}$. An argument similar to
(1) shows that
$$y\equiv(\prod_{i=1}^{m-\!1}\delta_{i}^{n(\bar{i},B,G)^{+}}\prod_{i=0}^{m-\!1}\delta_{i}^{n(\bar{i},B,G)^{-}}
\prod_{j=1}^{n-\!2\widetilde{k}}\xi_{i_{j}}^{k_{j}})B^{'}\otimes
G_{2}\otimes C_{1,1}^{I}$$ (mod
$\CTL(D)_{m,n}^{<(\widetilde{k},I)^{+}}$), where both $B^{'}$ and
the coefficient of $B^{'}\otimes G_{2}\otimes C_{1,1}^{I}$ are
independent of $G_{2}$.

(3) For each $a=B=B_{1}\otimes B_{2}\otimes\omega\in Q^{+}(n,k_{1})$
and each $C_{S,T}^{\lambda}=G=G_{1}\otimes G_{2}\otimes
C^{J}_{1,1}\in C(k_{2},J)^{-}$, it is always true that $y=B\cdot G
\in \CTL(D)_{m,n}^{<(k_{2},J)^{-}}$ since there is always a closed
cycle in $y$.

(4) For each $a=B=B_{1}\otimes B_{2}\otimes\omega\in Q^{-}(n,k_{1})$
and each $C_{S,T}^{\lambda}=G=G_{1}\otimes G_{2}\otimes
C^{J}_{1,1}\in C(k_{2},J)^{-}$, let $y=B\cdot G$. It is clear that
$y\in \CTL(D)_{m,n}^{<(k_{2},J)^{-}}$ when $k_{1}>k_{2}$. In the
case $k_{1}\leq k_{2}$, we also have $y\in
\CTL(D)_{m,n}^{<(k_{2},J)^{-}}$ if the horizontal arcs in $B_{2}$
and $G_{1}$ form an interior closed cycle with one blob. Otherwise
we have $y\in \CTL(D)_{m,n}^{\leqq(k_{2},J)^{-}}$. An argument
similar to (1) shows that
$$y\equiv(\prod_{i=1}^{m-\!1}\delta_{i}^{n(\bar{i},B,G)^{+}}\prod_{i=0}^{m-\!1}\delta_{i}^{n(\bar{i},B,G)^{-}}
\prod_{j=1}^{n-\!2k_{2}}\xi_{i_{j}}^{k_{j}})B^{'}\otimes
G_{2}\otimes C_{1,1}^{J}$$ (mod $\CTL(D)_{m,n}^{<(k_{2},J)^{-}}$),
where both $B^{'}$ and the coefficient of $B^{'}\otimes G_{2}\otimes
C_{1,1}^{J}$ are independent of $G_{2}$.

When $n$ is even, we have the following additional cases to prove.

(5) For each $a=B=B_{1}\otimes B_{2}\otimes\omega\in
Q^{-}_{i}(n,\frac{n}{2})\ (i=1,2)$ and each
$C_{S,T}^{\lambda}=G=G_{1}\otimes G_{2}\otimes C^{I}_{1,1}\in
C(\widetilde{k},I)^{+}$, let $y=B\cdot G$. It is clear that $y\in
\CTL(D)_{m,n}^{<\widetilde{k},I)^{+}}$ when
$\widetilde{k}<\frac{n}{2}$. In the case
$\widetilde{k}=\frac{n}{2}$, we have $$y\equiv
(\prod_{i=1}^{m-\!1}\delta_{i}^{n(\bar{i},B,G)^{+}}\prod_{i=0}^{m-\!1}\delta_{i}^{n(\bar{i},B,G)^{-}})
B^{'}_{1}\otimes G_{2}\otimes C_{1,1}^{\emptyset}\ \in
\CTL(D)_{m,n}^{(\widetilde{k},I)^{+}},$$ where $B^{'}_{1}\in
D^{+}(n,\widetilde{k})$ is determined by $B_{1}$ and the coefficient
of $B^{'}_{1}\otimes G_{2}\otimes C_{1,1}^{\emptyset}$ is completely
determined by the horizontal arcs in $B_{2}$ and $G_{1}$, thus both
$B^{'}_{1}$ and the coefficient of $B^{'}_{1}\otimes G_{2}\otimes
C_{1,1}^{\emptyset}$ are independent of $G_{2}$.

(6) For each $a=B=B_{1}\otimes B_{2}\otimes\omega\in
Q^{-}_{i}(n,\frac{n}{2})\ (i=1,2)$ and each
$C_{S,T}^{\lambda}=G=G_{1}\otimes G_{2}\otimes C^{J}_{1,1}\in
C(k,J)^{-}$, we always have $y=B\cdot G \in
\CTL(D)_{m,n}^{<(k,J)^{-}}$ since $k<\frac{n}{2}$.

(7) For each $a=B=B_{1}\otimes B_{2}\otimes\omega\in
Q^{-}_{1}(n,\frac{n}{2})$ (respectively, $Q^{-}_{2}(n,\frac{n}{2})$)
and each $C_{S,T}^{\lambda}=G=G_{1}\otimes G_{2}\otimes
C^{\emptyset}_{1,1}\in C(\frac{n}{2},\emptyset)^{-}_{1}$
(respectively, $C(\frac{n}{2},\emptyset)^{-}_{2}$), let $y=B\cdot
G$. If the horizontal arcs in $B_{2}$ and $G_{1}$ form an interior
closed cycle with one blob, then it is clear that $y \in
\CTL(D)_{m,n}^{<(\frac{n}{2},\emptyset)^{-}_{1}}$ (respectively,
$\CTL(D)_{m,n}^{<(\frac{n}{2},\emptyset)^{-}_{2}}$). Otherwise, we
have $$y=
(\prod_{i=1}^{m-\!1}\delta_{i}^{n(\bar{i},B,G)^{+}}\prod_{i=0}^{m-\!1}\delta_{i}^{n(\bar{i},B,G)^{-}})
B_{1}\otimes G_{2}\otimes C_{1,1}^{\emptyset}\ \in
\CTL(D)_{m,n}^{(\frac{n}{2},\emptyset)^{-}_{1}}$$ (respectively,
$\CTL(D)_{m,n}^{(\frac{n}{2},\emptyset)^{-}_{2}}$), where the
coefficient of $B_{1}\otimes G_{2}\otimes C_{1,1}^{\emptyset}$ is
completely determined by the horizontal arcs in $B_{2}$ and $G_{1}$,
thus is independent of $G_{2}$.

(8) For each $a=B=B_{1}\otimes B_{2}\otimes\omega\in
Q^{-}_{1}(n,\frac{n}{2})$ (respectively, $Q^{-}_{2}(n,\frac{n}{2})$)
and each $C_{S,T}^{\lambda}=G=G_{1}\otimes G_{2}\otimes
C^{\emptyset}_{1,1}\in C(\frac{n}{2},\emptyset)^{-}_{2}$
(respectively, $C(\frac{n}{2},\emptyset)^{-}_{1}$), it is always
true that $y=B\cdot G \in
\CTL(D)_{m,n}^{<(\frac{n}{2},\emptyset)^{-}_{2}}$ (respectively,
$\CTL(D)_{m,n}^{<(\frac{n}{2},\emptyset)^{-}_{1}}$) since there must
be an interior closed cycle with one blob in $y$, which is formed by
the horizontal arcs in $B_{2}$ and $G_{1}$.

(9) For each $a=B=B_{1}\otimes B_{2}\otimes\omega\in Q^{+}(n,k)$ and
each $C_{S,T}^{\lambda}=G=G_{1}\otimes G_{2}\otimes
C^{\emptyset}_{1,1}\in C(\frac{n}{2},\emptyset)^{-}_{1}$
(respectively, $C(\frac{n}{2},\emptyset)^{-}_{2}$), we always have
$y=B\cdot G \in \CTL(D)_{m,n}^{<(\frac{n}{2},\emptyset)^{-}_{1}}$
(respectively, $\CTL(D)_{m,n}^{<(\frac{n}{2},\emptyset)^{-}_{2}}$)
since there is always a closed cycle with one blob in $y$.

(10) For each $a=B=B_{1}\otimes B_{2}\otimes\omega\in Q^{-}(n,k)$
and each $C_{S,T}^{\lambda}=G=G_{1}\otimes G_{2}\otimes
C^{\emptyset}_{1,1}\in C(\frac{n}{2},\emptyset)^{-}_{1}$
(respectively, $C(\frac{n}{2},\emptyset)^{-}_{2}$), let $y=B\cdot
G$. If the horizontal arcs in $B_{2}$ and $G_{1}$ form an interior
closed cycle with one blob, then it is clear that $y \in
\CTL(D)_{m,n}^{<(\frac{n}{2},\emptyset)^{-}_{1}}$ (respectively,
$\CTL(D)_{m,n}^{<(\frac{n}{2},\emptyset)^{-}_{2}}$). Otherwise, we
have $$y= (\
\prod_{i=1}^{m-\!1}\delta_{i}^{n(\bar{i},B,G)^{+}}\prod_{i=0}^{m-\!1}\delta_{i}^{n(\bar{i},B,G)^{-}})
\ B^{'}\otimes G_{2}\otimes C_{1,1}^{\emptyset}\ \in \
\CTL(D)_{m,n}^{(\frac{n}{2},\emptyset)^{-}_{1}}$$ (respectively,
$\CTL(D)_{m,n}^{(\frac{n}{2},\emptyset)^{-}_{2}}$), where $B^{'}$ is
determined only by $B_{1}$, $B_{2}$ and $G_{1}$, and the coefficient
of $B^{'}\otimes G_{2}\otimes C_{1,1}^{\emptyset}$ is completely
determined by the horizontal arcs in $B_{2}$ and $G_{1}$. Therefore,
both $B^{'}$ and the coefficient of $B^{'}\otimes G_{2}\otimes
C_{1,1}^{\emptyset}$ are independent of $G_{2}$.

The above argument implies that the condition $(C3)$ in Definition
\ref{D4.1.1} follows. Thus the proof of the theorem is completed.\
\hfill $\square$

\bigskip

{\it Remark}. In particular, we obtain from the above theorem that
Temperley-Lieb algebras of type $D$ over an arbitrary field are
cellular.

\section{Irreducible representations and
quasi-heredity of $\CTL(D)_{m,n}$}

$\indent$The theory of cellular algebras can help us to determine
all the irreducible representations of a cellular algebra. Cellular
algebras and quasi-hereditary algebras are closely related. The
purpose of this section is to investigate the irreducible
representations of cyclotomic Temperley-Lieb algebras of type $D$
and then determine which parameters yield a quasi-hereditary
cyclotomic Temperley-Lieb algebra of type $D$.

We first recall that given a cellular algebra $A$ with cell datum
$(\Lambda,M,C,i)$, one can define for each $\lambda\in\Lambda$ a
\emph{cell module} $W(\lambda)$ and a bilinear form
$\phi_{\lambda}:\ W(\lambda)\otimes_{K}W(\lambda)\rightarrow K$ as
follows. As a vector space, $W(\lambda)$ has a $K$-basis
$\{C_S^{\lambda}\,|\,S\in M(\lambda)\}$, and the module structure is
given by
$$a\,C_S^{\lambda}=\sum_{U\in M(\lambda)}r_a(U,S)C^{\lambda}_{U}\
,$$where the coefficients $r_a(U,S)$ are determined by $(C3)$ in
Definition \ref{D4.1.1}. The bilinear form $\phi_{\lambda}$ is
defined by
$$\phi_{\lambda}(C_S^{\lambda},C_T^{\lambda})C_{U,V}^{\lambda}
\equiv C_{U,S}^{\lambda}C_{T,V}^{\lambda}$$modulo the ideal
generated by all basis elements with upper index strictly smaller
than $\lambda$.

The theory of cellular algebras (see Graham and Lehrer, 1996) shows
that the isomorphism classes of simple $A$-modules are parametrized
by the set
$\Lambda_0=\{\lambda\in\Lambda\,|\,\phi_{\lambda}\neq0\}$. It can be
realized in the following way. We write rad$(\lambda)$ for the
subspace of the cell module $W(\lambda)$ given by $\{x\in
W(\lambda)\,|\,\phi_{\lambda}(x,y)=0\ \mbox{for all}\ y\in
W(\lambda)\}$. It has been proven that rad($\lambda$) is a submodule
of $W(\lambda)$. As a result, the factor module
$W(\lambda)/\mbox{rad}(\lambda)$ with $\lambda\in\Lambda_0$ gives
rise to a simple $A$-module. In this case, we write $S(\lambda)$ for
$W(\lambda)/\mbox{rad}(\lambda)$.

The following theorem is a parametrization of the irreducible
representations of cyclotomic Temperley-Lieb algebras of type $D$.
To state the theorem, we introduce the following notations.
$$S_{1}=\left\{S((\widetilde{k},I)^{+}),\ S((k,J)^{-})\,\left|\,
\begin{array}{c}
1\leq \widetilde{k}\leq[\frac{n}{2}],\ 0\leq k\leq[\frac{n}{2}],
\\I=(i_{1},\cdots,i_{n-\!2\widetilde{k}})\in \Lambda(m,n-\!2\widetilde{k}),
\\J=(j_{1},\cdots,j_{n-\!2k})\in \Lambda(m,n-\!2k)
\\ with\ all\ i_{h},\ j_{h}\ divisible\ by\ p^{t}
\end{array}
\right.\right\}$$
$$S_{2}=\left\{S((\widetilde{k},I)^{+}),\ S((k,J)^{-})\,\left|\,
\begin{array}{c}
1\leq \widetilde{k}\leq[\frac{n}{2}]-1,\ 0\leq k\leq[\frac{n}{2}]-1,
\\I=(i_{1},\cdots,i_{n-\!2\widetilde{k}})\in \Lambda(m,n-\!2\widetilde{k}),
\\J=(j_{1},\cdots,j_{n-\!2k})\in \Lambda(m,n-\!2k)
\\ with\ all\ i_{h},\ j_{h}\ divisible\ by\ p^{t}
\end{array}
\right.\right\}$$\\
$$S_{2}^{'}=\biggl\{S((\frac{n}{2},\emptyset)^{+}),\
S((\frac{n}{2},\emptyset)_{1}^{-}),\
S((\frac{n}{2},\emptyset)_{2}^{-})\biggr\}$$

\begin{Thm}
Let $K$ be a splitting field of $x^{m}-1$ and $\Char K=p$. Write
$m=p^{t}s$ with $(p,s)=1$ and $t\geq0$ (in the case $p=0$, set
$0^{0}=1$). Then we have the following.

(1) If $n$ is odd, then the set $S_{1}$ forms a complete set of
non-isomorphic simple $\CTL(D)_{m,n}$-modules.

(2) If $n$ is even.

\ \ ($i$) If not all $\delta_{i}$ are zero, then the set $S_{2}\cup
S_{2}^{'}$ forms a complete set of non-isomorphic simple
$\CTL(D)_{m,n}$-modules.

\ \ ($ii$) If all $\delta_{i}$ are zero, then the set $S_{2}$ is a
complete set of non-isomorphic simple
$\CTL(D)_{m,n}$-modules.\label{T5.1.1}
\end{Thm}

To prove Theorem \ref{T5.1.1}, we first recall the following lemma
from Rui and Xi (2004) which describes the simple
$\G_{m,n}$-modules.

\begin{Lem}{\rm (Rui and Xi, 2004)}. Let $K$ be a splitting field of $x^{m}-1$ and $\Char K=p$.

(1) If $p$ divides $m$, say, $m=p^{t}s$ with $(p,s)=1$, then
$\{S(I)\,|\,I=(i_{1},\cdots,i_{n})$ with all $i_{j}$ divisible by
$p^{t}\}$ forms a complete set of non-isomorphic simple
$\G_{m,n}$-modules whose cardinality is $s^{n}$.

(2) If $p$ does not divide $m$, then $\{S(I)\,|\,I\in\Lambda(m,n)\}$
is a complete set of non-isomorphic simple $\G_{m,n}$-modules. In
this case, the algebra $\G_{m,n}$ is semisimple.\label{L5.2.1}
\end{Lem}

For any $I\in \Lambda(m,n-\!2k)$, let $\psi_{I}$ be the bilinear
form defined on the cell module of $G_{m,n-\!2k}$ associated with
the index $I$. Similarly, for any $\lambda\in \Lambda_{m,n}$, let
$\phi_{\lambda}$ be the bilinear form defined on the cell module
$W(\lambda)$ of $\CTL(D)_{m,n}$. To study the irreducible
representations of $\CTL(D)_{m,n}$, we first discuss when the
bilinear form $\phi_{\lambda}$ is equal to zero.

\begin{Lem}(1) If $\psi_{I}=0$ for some $I\in \Lambda(m,n-\!2k)$ and
$1\leq k\leq [n/2]$, $k\neq\frac{n}{2}$, then $\phi_{(k,I)^{+}}=0$.

(2) If $\psi_{I}=0$ for some $I\in \Lambda(m,n-\!2k)$ and $0\leq
k\leq [n/2]$, $k\neq\frac{n}{2}$, then
$\phi_{(k,I)^{-}}=0$.\label{L5.2.2}
\end{Lem}
{\it Proof.} For any $v_{1},v_{2}\in D^{+}(n,k)$ and $I\in
\Lambda(m,n-\!2k)$, let $y=(v_{1}\otimes v_{2}\otimes
C_{1,1}^{I})\cdot(v_{1}\otimes v_{2}\otimes C_{1,1}^{I})$. If $y\in
\CTL(D)_{m,n}^{<(k,I)^{+}}$, then $\phi_{(k,I)^{+}}=0$. Otherwise,
suppose $y\equiv v_{1}\otimes v_{2}\otimes x C_{1,1}^{I}C_{1,1}^{I}$
(mod $\CTL(D)_{m,n}^{<(k,I)^{+}}$), where $x\in G_{m,n-\!2k}$. If
$\psi_{I}=0$, then $C_{1,1}^{I}C_{1,1}^{I}\equiv0$ (mod
$G_{m,n-\!2k}^{\prec I}$). Hence $x C_{1,1}^{I}C_{1,1}^{I}\equiv0$
(mod $G_{m,n-\!2k}^{\prec I}$) and $v_{1}\otimes v_{2}\otimes x
C_{1,1}^{I}C_{1,1}^{I}\equiv0$ (mod $\CTL(D){m,n}^{<(k,I)^{+}}$),
thus $\phi_{(k,I)^{+}}=0$. By a similar argument, we have
$\phi_{(k,I)^{-}}=0$ as required.\ \hfill $\square$

\begin{Lem}Assume that $\psi_{I}\neq0$ for some $I\in \Lambda(m,n-\!2k)$. Then
$\phi_{(k,I)^{+}}\neq0$ and $\phi_{(k,I)^{-}}\neq0$ if $n$ is odd or
$n$ is even but $k\neq \frac{n}{2}$. \label{L5.2.3}
\end{Lem}
{\it Proof.} We can show that $\phi_{(k,I)^{+}}$ and
$\phi_{(k,I)^{-}}$ are not zero by choosing some special
$m$-cyclotomic $D$-admissible dangles of type $(n,k)$. Let $D_{1}$
be the $m$-cyclotomic $D$-admissible dangle of type I (respectively,
type II) with horizontal arcs $\{1<2\},\
\{3<4\},\cdots,\{2k-\!1<2k\}$, and let $D_{2}$ be the $m$-cyclotomic
$D$-admissible dangle of type I (respectively, type II) with
horizontal arcs $\{2<3\},\ \{4<5\},\cdots,\{2k<2k+\!1\}$. Each
horizontal arc does not carry dots and there are no blobs on each
horizontal arc and vertical line. It is clear that $(D_{1}\otimes
D_{2}\otimes C_{1,1}^{I})\cdot(D_{1}\otimes D_{2}\otimes
C_{1,1}^{I})=D_{1}\otimes D_{2}\otimes C_{1,1}^{I}C_{1,1}^{I}$, thus
we get $\phi_{(k,I)^{+}}\neq0$ (respectively,
$\phi_{(k,I)^{-}}\neq0$) from $\psi_{I}\neq0$.\ \hfill $\square$

\begin{Lem}Assume that $n$ is even and $k=\frac{n}{2}$. Then

(1) $\phi_{(\frac{n}{2},\emptyset)^{+}}=0$ if and only if all
$\delta_{i}$ are zero.

(2) $\phi_{(\frac{n}{2},\emptyset)^{-}_{i}}=0$ ($i=1,2$) if and only
if all $\delta_{i}$ are zero. \label{L5.2.4}
\end{Lem}
{\it Proof.} We only prove statement (2). The proof of statement (1)
is similar. We again choose some special $m$-cyclotomic
$D$-admissible dangles of type $(n,\frac{n}{2})$ to verify the
lemma. Let $D_{1}$ be the $m$-cyclotomic $D$-admissible of type II
with horizontal arcs $\{1<n\},\ \{2<3\},\
\{4<5\},\cdots,\{n-\!2<n-\!1\}$, and let $D_{2}$ be the
$m$-cyclotomic $D$-admissible dangle of type II with horizontal arcs
$\{1<2\},\ \{3<4\},\cdots,\{n-\!1<n\}$.

Suppose that $\delta_{i}\neq0$ for some $0\leq i\leq m-1$. Let
$D_{1}^{+}=D_{1}$ and $D_{2}^{+}$ to be the dangle decorated by $i$
dots in the horizontal arc $\{1<2\}$ in $D_{2}$. In this case, we
get $(D_{1}^{+}\otimes D_{2}^{+}\otimes
C_{1,1}^{\emptyset})\cdot(D_{1}^{+}\otimes D_{2}^{+}\otimes
C_{1,1}^{\emptyset})=\delta_{i}(D_{1}^{+}\otimes D_{2}^{+}\otimes
C_{1,1}^{\emptyset})$. So
$\phi_{(\frac{n}{2},\emptyset)^{-}_{1}}\neq0$. We take $D^{+}_{1}$
to be the dangle decorated by one blob in the horizontal arc
$\{1<n\}$ in $D_{1}$ and $D_{2}^{+}$ to be the dangle decorated by
one blob and $i$ dots in the horizontal arc $\{1<2\}$ in $D_{2}$. In
this case, we also get $(D_{1}^{+}\otimes D_{2}^{+}\otimes
C_{1,1}^{\emptyset})\cdot(D_{1}^{+}\otimes D_{2}^{+}\otimes
C_{1,1}^{\emptyset})=\delta_{i}(D_{1}^{+}\otimes D_{2}^{+}\otimes
C_{1,1}^{\emptyset})$, so
$\phi_{(\frac{n}{2},\emptyset)^{-}_{2}}\neq0$.

Suppose that $\delta_{i}=0$ for all $0\leq i\leq m-1$. For any
$v_{1},v_{2}\in D^{-}_{i}(n,\frac{n}{2})$ ($i=1,2$), the composition
of $v_{1}\otimes v_{2}\otimes C_{1,1}^{\emptyset}$ with itself must
contain at least one interior closed cycle since there are no
vertical lines in the dangles $v_{1}$ and $v_{2}$. Note that these
interior closed cycles will provide some zero factors $\delta_{i}$.
Thereby the product $(v_{1}\otimes v_{2}\otimes
C_{1,1}^{\emptyset})\cdot(v_{1}\otimes v_{2}\otimes
C_{1,1}^{\emptyset})$ is zero, and thus
$\phi_{(\frac{n}{2},\emptyset)^{-}_{i}}=0$ ($i=1,2$).\ \hfill
$\square$

\bigskip

Now we are in the position to prove Theorem \ref{T5.1.1}.

{\it Proof of Theorem \ref{T5.1.1}.} Let
$(\widetilde{k},I)^{+}\in\Lambda_{m,n}$ and
$(k,J)^{-}\in\Lambda_{m,n}$. In the case
$(\widetilde{k},I)^{+}=(\frac{n}{2},\emptyset)^{+}$ (respectively,
$(k,J)^{-}=(\frac{n}{2},\emptyset)_{i}^{-}$ $(i=1,2)$), we know by
Lemma \ref{L5.2.4} that $\phi_{(\frac{n}{2},\emptyset)^{+}}=0$
(respectively, $\phi_{(\frac{n}{2},\emptyset)_{i}^{-}}=0$ $(i=1,2)$)
if and only if $\delta_{i}=0$ for all $0\leq i\leq m-1$. If
$(\widetilde{k},I)^{+}\neq(\frac{n}{2},\emptyset)^{+}$
(respectively, $(k,J)^{-}\neq(\frac{n}{2},\emptyset)_{i}^{-}$
$(i=1,2)$), then by Lemma \ref{L5.2.2} and Lemma \ref{L5.2.3} we see
that $\phi_{(\widetilde{k},I)^{+}}=0$ (respectively,
$\phi_{(k,J)^{-}}=0$) if and only if $\psi_{I}=0$ (respectively,
$\psi_{J}=0$). Moreover, by Lemma \ref{L5.2.1} we know that
$\psi_{I}\neq0$ if and only if $p^{t}$ divides all $i_{h}$, where
$i_{h}$ are the components of
$I=(i_{1},i_{2},\cdots,i_{n-\!2\widetilde{k}})\in\Lambda(m,n-\!2\widetilde{k})$
and $\psi_{J}\neq0$ if and only if $p^{t}$ divides all $j_{h}$,
where $j_{h}$ are the components of
$J=(j_{1},j_{2},\cdots,j_{n-\!2k})\in\Lambda(m,n-\!2k)$. This
finishes the proof.\ \hfill $\square$

\bigskip

In the rest of this section we will determine which parameters yield
a quasi-hereditary cyclotomic Temperley-Lieb algebra of type $D$.
Quasi-hereditary algebras are used to describe the highest weight
categories appearing in the representation theory of semisimple Lie
algebras and algebraic groups. We first recall the definition of
quasi-hereditary algebras.

\begin{Def}{\rm (Cline, Parshall and Scott, 1988). Let $A$ be a
 $K$-algebra. An ideal $J$ of $A$ is called a \emph{heredity ideal}
 if $J$ is idempotent, $J(radA)J=0$, and $J$ is a projective left
 $($or, right$)$ $A$-module. The algebra $A$ is called
 \emph{quasi-hereditary} provided there is a finite chain $0=J_{0}\subset
J_{1}\subset \cdots \subset J_{n}=A$ of ideals in $A$ such that
$J_{j}/J_{j-\!1}$ is a heredity ideal of $A/J_{j-\!1}$ for all $j$.
Such a chain is then called a \emph{heredity chain} of the
quasi-hereditary algebra $A$.}\label{D6.1.1}
\end{Def}

As indicated by the ideals appearing in cell chains, there are close
connections between cellular algebras and quasi-hereditary algebras.
The class of cellular algebras has a large intersection with the
class of quasi-hereditary algebras. Typical examples of
quasi-hereditary algebras obtained from cellular algebras include
Temperley-Lieb algebras of type $A$ with non-zero parameters (see
Westbury, 1995) and Birman-Wenzl algebras for most choices of
parameters (see Xi, 2000). Recently, Xi (2002) proved that a
cellular algebra is quasi-hereditary if and only if the first
cohomology groups between cell modules and dual modules are always
trivial. A stronger statement describing the quasi-heredity of
cellular algebras was given later by Cao (2003). One can refer to
the survey paper by K\"{o}nig and Xi (1999c) for a comparison
between cellular algebras and quasi-hereditary algebras. Now we
recall the following theorem which determines those cellular
algebras which are quasi-hereditary.

\begin{Thm}{\rm (K\"onig and Xi, 1999a)}. Let $K$ be a field and $A$ a cellular $K$-algebra
(with respect to an involution $i$). Then the following statements
are equivalent.

(1) Some cell chain of $A$ (with respect to some involution,
possibly different from $i$) is a heredity chain as well, thus it
makes $A$ into a quasi-hereditary algebra.

(1') There is a cell chain of $A$ (with respect to some involution,
possibly different from $i$) whose length equals the number of
isomorphism classes of simple $A$-modules.

(2) $A$ has finite global dimension.

(3) The Cartan matrix of $A$ has determinant one.

(4) Any cell chain of $A$ (with respect to any involution) is a
heredity chain.\label{T6.1.2}
\end{Thm}

As an immediate consequence of the above theorem, we have the
following corollary.

\begin{Coro} Let $K$ be a field and $A$ a cellular $K$-algebra.
Then $A$ is quasi-hereditary if and only if
$|\Lambda|=|\Lambda_{0}|$. \label{C6.1.3}
\end{Coro}
{\it Proof.} Suppose that $|\Lambda|\neq|\Lambda_{0}|$. Then there
exists a cell chain which is not a heredity chain since there exists
a cell ideal which is not a heredity one. So $A$ is not
quasi-hereditary from the statement (4) in Theorem \ref{T6.1.2}.

If $|\Lambda|=|\Lambda_{0}|$, then there exists a cell chain whose
length equals the number of isomorphism classes of simple
$A$-modules. So $A$ is quasi-hereditary from the statement (1') in
Theorem \ref{T6.1.2}.\ \hfill $\square$

\bigskip

The following theorem gives a necessary and sufficient condition for
a cyclotomic Temperley-Lieb algebras of type $D$ to be
quasi-hereditary.

\begin{Thm}
Let $K$ be a splitting field of $x^{m}-1$ and let $\Char K=p$. Then
the cyclotomic Temperley-Lieb algebra of type $D$ is
quasi-hereditary if and only if $p$ does not divide $m$ and one of
the following holds:

(1) $n$ is odd.

(2) $n$ is even, but not all $\delta_{i}$ are zero.\label{T6.2.1}
\end{Thm}
{\it Proof.} By Corollary \ref{C6.1.3}, we know that $\CTL(D)_{m,n}$
is quasi-hereditary if and only if the index set of cell modules of
$\CTL(D)_{m,n}$ coincides with that of simple modules. Moreover, the
coincidence occurs if and only if $p$ does not divide $m$, and
either $n$ is odd or not all $\delta_{i}$ are zero by Theorem
\ref{T5.1.1}.\ \hfill $\square$

\bigskip

For the case not displayed in Theorem \ref{T6.2.1}, we may get a
quasi-hereditary quotient of $\CTL(D)_{m,n}$.

\begin{Prop}
Suppose that $K$ is a splitting field of $x^{m}-1$ and that $p$ does
not divide $m$. If $n$ is even and all $\delta_{i}$ are zero, then
the factor algebra $\CTL(D)_{m,n}/J$ is quasi-hereditary, where $J$
is the ideal of $\CTL(D)_{m,n}$ generated by all $m$-cyclotomic
$D$-admissible diagrams without vertical arcs.
\end{Prop}
{\it Proof.} If $n$ is even and all $\delta_{i}$ are zero, then the
ideal $J$ generated by all $m$-cyclotomic $D$-admissible diagrams
without vertical arcs is nilpotent by Lemma \ref{L5.2.4}. It is
clear that the factor algebra $\CTL(D)_{m,n}/J$ is again cellular
with respect to the induced involution. In addition, under the
assumption of the proposition the index set of cell modules of
$\CTL(D)_{m,n}/J$ coincides with that of simple
$\CTL(D)_{m,n}/J$-modules. Therefore $\CTL(D)_{m,n}/J$ is a
quasi-hereditary algebra.\ \hfill $\square$

\section{An example}

$\indent$ In this section, we will go though a concrete example of
cyclotomic Temperley-Lieb algebras of type $D$ in the case $m=2$ and
$n=4$ to illustrate the results in the previous sections.

In the case $n=4$, the $K$-dimension of Temperley-Lieb algebra of
type $D$ is $48$ with $13$ $D$-admissible diagrams of type I and
$35$ $D$-admissible diagrams of type II. In the case $m=2$ and
$n=4$, the number of $2$-cyclotomic $D$-admissible diagrams of type
I is $208$ and the number of $m$-cyclotomic $D$-admissible diagrams
of type II is $560$. This is a total of $768$ which is the
$K$-dimension of $\CTL(D)_{2,4}$.

We first describe the cellular structure of $\CTL(D)_{2,4}$. Recall
that $\Lambda(2,4)=\{(j_{1},j_{2},j_{3},j_{4})\,|\,1\leq j_{h}\leq
2\}$ and $\Lambda(2,2)=\{(i_{1},i_{2})\,|\,1\leq i_{h}\leq 2\}$. The
partial orderings of $(\Lambda(2,4),\preceq)$ and
$(\Lambda(2,2),\preceq)$ are illustrated by the following figures.
Two elements connected in the figures are comparable and the larger
element sits above the smaller.
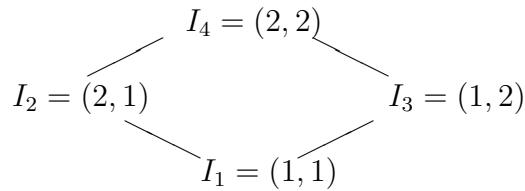
\begin{center}
\setlength{\unitlength}{1cm}
\begin{picture}(10,9.5)
\put(4.5,0.5){$(1,1,1,1)$} \put(5,0){$J_{1}$}
\put(0,2.0){$(1,1,1,2)$} \put(2.7,2.0){$(1,1,2,1)$}
\put(5.4,2.0){$(1,2,1,1)$} \put(8.2,2.0){$(2,1,1,1)$}
\put(0.8,1.5){$J_{2}$} \put(3.5,1.5){$J_{3}$} \put(6.2,1.5){$J_{4}$}
\put(9,1.5){$J_{5}$} \put(4.5,0.8){\line(-2,1){1}}
\put(4.5,0.8){\line(-4,1){3}} \put(5.8,0.8){\line(2,1){1}}
\put(5.8,0.8){\line(4,1){3}} \put(0,4.5){$(2,2,1,1)$}
\put(0.8,4){$J_{6}$} \put(2,4.5){$(2,1,2,1)$} \put(2.8,4){$J_{7}$}
\put(4,4.5){$(2,1,1,2)$} \put(4.8,4){$J_{8}$}
\put(6,4.5){$(1,2,2,1)$} \put(6.8,4){$J_{9}$}
\put(8,4.5){$(1,2,1,2)$} \put(8.8,4){$J_{10}$}
\put(10,4.5){$(1,1,2,2)$} \put(10.8,4){$J_{11}$}

\put(1.5,2.3){\line(5,2){3.5}} \put(1.5,2.3){\line(5,1){7.05}}
\put(1.5,2.3){\line(6,1){9}} \put(3.5,2.3){\line(5,2){3.4}}
\put(3.5,2.3){\line(-2,5){0.55}}
 \put(3.5,2.3){\line(5,1){7}}
\put(6.5,2.3){\line(3,2){2.05}} \put(6.5,2.3){\line(1,3){0.45}}

\put(6.5,2.3){\line(-4,1){5.5}} \put(8.5,2.3){\line(-5,2){3.5}}
\put(8.5,2.3){\line(-4,1){5.5}} \put(8.5,2.3){\line(-5,1){7.5}}
\put(0,7){$(2,2,2,1)$} \put(2.8,7){$(2,2,1,2)$}
\put(5.6,7){$(2,1,2,2)$} \put(8.4,7){$(1,2,2,2)$}
\put(0.8,6.5){$J_{12}$} \put(3.6,6.5){$J_{13}$}
\put(6.4,6.5){$J_{14}$} \put(9.2,6.5){$J_{15}$}
\put(1.5,6.2){\line(-1,-2){0.65}} \put(1.5,6.2){\line(1,-1){1.25}}
 \put(1.5,6.2){\line(4,-1){5.05}}

 \put(3.5,6.2){\line(-2,-1){2.6}}
\put(3.5,6.2){\line(4,-3){1.7}} \put(3.5,6.2){\line(4,-1){5}}
 \put(6.5,6.2){\line(-1,-1){1.3}}
\put(6.5,6.2){\line(-3,-1){3.75}}  \put(6.5,6.2){\line(3,-1){4}}
\put(8.5,6.2){\line(0,-1){1.3}} \put(8.5,6.2){\line(3,-2){2}}
\put(8.5,6.2){\line(-3,-2){1.95}} \put(4.5,8.7){$(2,2,2,2)$}
\put(5,8.2){$J_{16}$} \put(4.5,8){\line(-4,-1){2.8}}
\put(4.5,8){\line(-2,-1){1}} \put(5.8,8){\line(2,-1){1}}
\put(5.8,8){\line(5,-1){3}}
\end{picture}
\end{center}
\begin{center}
Figure 6. Partially ordered set $(\Lambda(2,4),\preceq)$
\end{center}
\begin{center}
\setlength{\unitlength}{1cm}
\begin{picture}(10,4)
\put(4.5,1){$I_{1}=(1,1)$} \put(2,2){$I_{2}=(2,1)$}
\put(7,2){$I_{3}=(1,2)$} \put(4.3,3){$I_{4}=(2,2)$}
\put(4.5,1.3){\line(-2,1){1}} \put(5.8,1.3){\line(2,1){1}}
\put(3,2.4){\line(2,1){1}} \put(7,2.4){\line(-2,1){1}}
\end{picture}
\end{center}
\begin{center}
Figure 7. Partially ordered set $(\Lambda(2,2),\preceq)$
\end{center}

The partial ordering of $(\Lambda_{2,4},\leqq)$ can be illustrated
by the following figure.
\begin{center}
\setlength{\unitlength}{1cm}
\begin{picture}(10,15)
\put(4.5,0){$(2,\emptyset)^{+}$} \put(5,0.3){\line(0,1){0.5}}
\put(4.5,1){$(1,I_{1})^{+}$} \put(2,2){$(1,I_{2})^{+}$}
\put(7,2){$(1,I_{3})^{+}$} \put(4.5,3){$(1,I_{4})^{+}$}
\put(4.5,1.3){\line(-2,1){1.3}} \put(5.8,1.3){\line(2,1){1.2}}
\put(3,2.4){\line(2,1){1.3}} \put(7,2.4){\line(-2,1){1.3}}
\put(2,4){$(2,\emptyset)^{-}_{1}$}
\put(7,4){$(2,\emptyset)^{-}_{2}$} \put(4.5,3.3){\line(-2,1){1.3}}
\put(5.8,3.3){\line(2,1){1.2}} \put(3,4.4){\line(2,1){1.3}}
\put(7,4.4){\line(-2,1){1.3}} \put(4.5,5){$(1,I_{1})^{-}$}
\put(2,6){$(1,I_{2})^{-}$} \put(7,6){$(1,I_{3})^{-}$}
\put(4.5,7){$(1,I_{4})^{-}$} \put(4.5,5.3){\line(-2,1){1.3}}
\put(5.8,5.3){\line(2,1){1.2}} \put(3,6.4){\line(2,1){1.3}}
\put(7,6.4){\line(-2,1){1.3}} \put(4.5,8){$(0,J_{1})^{-}$}
\put(5,7.3){\line(0,1){0.5}} \put(1,9){$(0,J_{2})^{-}$}
\put(3,9){$(0,J_{3})^{-}$} \put(6,9){$(0,J_{4})^{-}$}
\put(8,9){$(0,J_{5})^{-}$} \put(4.5,8.3){\line(-2,1){1}}
\put(4.5,8.3){\line(-4,1){2.5}} \put(5.8,8.3){\line(2,1){1}}
\put(5.8,8.3){\line(4,1){2.5}} \put(0,11){$(0,J_{6})^{-}$}
\put(2,11){$(0,J_{7})^{-}$} \put(4,11){$(0,J_{8})^{-}$}
\put(6,11){$(0,J_{9})^{-}$} \put(8,11){$(0,J_{10})^{-}$}
\put(10,11){$(0,J_{11})^{-}$} \put(1.5,9.3){\line(5,2){3.5}}
 \put(1.5,9.3){\line(5,1){7.05}}
\put(1.5,9.3){\line(6,1){9}} \put(3.5,9.3){\line(5,2){3.4}}
 \put(3.5,9.3){\line(-2,5){0.55}}
\put(3.5,9.3){\line(5,1){7}}
 \put(6.5,9.3){\line(3,2){2.05}}

\put(6.5,9.3){\line(-4,1){5.5}} \put(6.5,9.3){\line(1,3){0.45}}

 \put(8.5,9.3){\line(-5,2){3.5}}
\put(8.5,9.3){\line(-4,1){5.5}} \put(8.5,9.3){\line(-5,1){7.5}}
\put(1,13){$(0,J_{12})^{-}$} \put(3,13){$(0,J_{13})^{-}$}
\put(6,13){$(0,J_{14})^{-}$} \put(8,13){$(0,J_{15})^{-}$}
\put(1.5,12.7){\line(-1,-2){0.65}} \put(1.5,12.7){\line(1,-1){1.25}}
 \put(1.5,12.7){\line(4,-1){5.05}}

\put(3.5,12.7){\line(-2,-1){2.6}} \put(3.5,12.7){\line(4,-3){1.7}}
 \put(3.5,12.7){\line(4,-1){5}}

\put(6.5,12.7){\line(-3,-1){3.75}} \put(6.5,12.7){\line(-1,-1){1.3}}
\put(6.5,12.7){\line(3,-1){4}} \put(8.5,12.7){\line(0,-1){1.3}}
\put(8.5,12.7){\line(3,-2){2}} \put(8.5,12.7){\line(-3,-2){1.95}}

 \put(4.5,14){$(0,J_{16})^{-}$}
\put(4.5,14){\line(-4,-1){2}} \put(4.5,14){\line(-2,-1){1}}
\put(5.8,14){\line(2,-1){1}} \put(5.8,14){\line(4,-1){2.5}}
\end{picture}
\end{center}
\begin{center}
Figure 8. Partially ordered set $(\Lambda_{2,4},\leqq)$
\end{center}
Suppose $v_{1},v_{2}\in D^{-}(4,1)$ are $2$-decorated dangles as
follows.
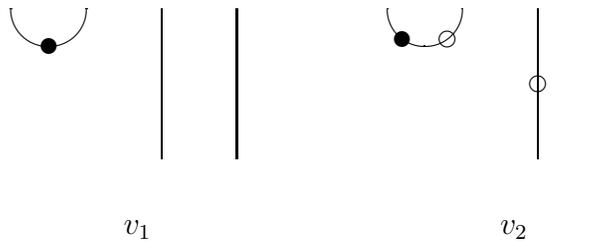
\begin{center}
\setlength{\unitlength}{1cm}
\begin{picture}(8,3)
\put(1.5,0){$v_{1}$} \put(6.5,0){$v_{2}$} \put(2,1){\line(0,1){2}}
\put(3,1){\line(0,1){2}} \put(7,1){\line(0,1){2}}
\put(8,1){\line(0,1){2}} \put(0.5,3){\oval(1,1)[b]}
\put(5.5,3){\oval(1,1)[b]} \put(0.5,2.5){\circle*{.2}}
\put(5.2,2.6){\circle*{.2}} \put(5.8,2.6){\circle{.2}}
\put(7,2){\circle{.2}}
\end{picture}
\end{center}
\begin{center}
Figure 9. $2$-Decorated dangles
\end{center}
For $(1,I_{2})^{-}\in \Lambda_{2,4}$ and
$(v_{1},I_{2}),(v_{2},I_{2})\in M((1,I_{2})^{-})$, we have
$C^{(1,I_{2})^{-}}_{(v_{1},I_{2}),(v_{2},I_{2})}=v_{1}\otimes
v_{2}\otimes C^{I_{2}}_{1,1}$ as follows.

\begin{center}
\setlength{\unitlength}{1cm}
\begin{picture}(8,3)
\put(3,0){$C^{(1,I_{2})^{-}}_{(v_{1},I_{2}),(v_{2},I_{2})}$}
\put(4,2){$+$} \put(2,1){\line(0,1){2}} \put(3,1){\line(0,1){2}}
\put(7,1){\line(0,1){2}} \put(8,1){\line(0,1){2}}
\put(0.5,3){\oval(1,1)[b]} \put(5.5,3){\oval(1,1)[b]}
\put(0.5,1){\oval(1,1)[t]} \put(5.5,1){\oval(1,1)[t]}
\put(0.5,2.5){\circle*{.2}} \put(5.5,2.5){\circle*{.2}}
\put(5.8,1.4){\circle{.2}} \put(7,2){\circle{.2}}
\put(0.2,1.4){\circle*{.2}} \put(0.8,1.4){\circle{.2}}
\put(5.2,1.4){\circle*{.2}} \put(2,2){\circle{.2}}
\put(3,2){\circle*{.2}}
\end{picture}
\end{center}
\begin{center}
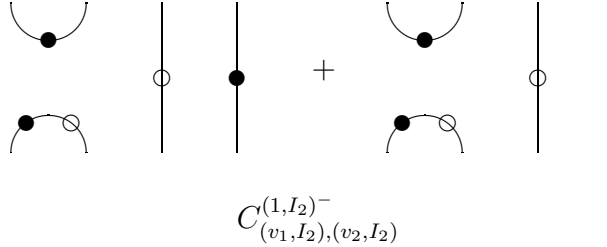
Figure 10. An element of the cellular basis
\end{center}
We can write all other elements of the cellular basis of
$\CTL(D)_{2,4}$ in a similar way.

Secondly, by Theorem \ref{T5.1.1} we can determine all the
irreducible representations of $\CTL(D)_{2,4}$. There are four
cases.

(1) If Char$K=0$ or Char$K=p\geq 3$ and $\delta_{0}\neq
0,\delta_{1}=0$ or $\delta_{0}=1,\delta_{1}\neq 0$, then the set
$\{S(\lambda)\,|\,\lambda\in\Lambda_{2,4}\}$ forms a complete set of
non-isomorphic simple $\CTL(D)_{m,n}$-modules.

(2) If Char$K=0$ or Char$K=p\geq 3$ and $\delta_{0}=\delta_{1}=0$,
then the set $\{S(\lambda)\,|\,\lambda\in\Lambda_{2,4}\}\backslash
\{S((2,\emptyset)^{+}),\ S((2,\emptyset)_{1}^{-}),\
S((2,\emptyset)_{2}^{-})\} $ forms a complete set of non-isomorphic
simple $\CTL(D)_{m,n}$-modules.

(3) If Char$K=2$ and $\delta_{0}\neq 0,\delta_{1}=0$ or
$\delta_{0}=1,\delta_{1}\neq 0$, then the set $\{S((1,I_{4})^{+}),$
$S((1,I_{4})^{-}),\ S((0,J_{16})^{-}), S((2,\emptyset)^{+}),\
S((2,\emptyset)_{1}^{-}),\ S((2,\emptyset)_{2}^{-})\}$ forms a
complete set of non-isomorphic simple $\CTL(D)_{m,n}$-modules.

(4) If Char$K=2$ and $\delta_{0}=\delta_{1}=0$, then the set
$\{S((1,I_{4})^{+}),\ S((1,I_{4})^{-}),$ $S((0,J_{16})^{-})\}$ forms
a complete set of non-isomorphic simple $\CTL(D)_{m,n}$-modules.

Finally, we know by Theorem \ref{T6.2.1} that $\CTL(D)_{2,4}$ is
quasi-hereditary if and only if Char$K=0$ or Char$K=p\geq 3$ and
$\delta_{0}\neq 0,\delta_{1}=0$ or $\delta_{0}=1,\delta_{1}\neq 0$.

 \bigskip \bigskip \bigskip

 {\Large {\bf Acknowledgements}}
  \bigskip

 The author would like to thank the referee
for his/her useful comments. The author is also grateful to
Professor Changchang Xi for his encouragement and advice.

 \bigskip \bigskip
 {\Large {\bf References}}
 \bigskip
\begin{enumerate}
\thispagestyle{plain}
 \footnotesize

\item Ariki, S., Koike, K. (1994). A Hecke algebra of $(\mathbb{Z}/r\mathbb{Z})\wr\frak{S}_n$
and construction of its irreducible representations. {\it Adv.
Math.} 106:216--243.\label{AK94}

\item Cao, Y. Z. (2003). On the quasi-heredity and the semi-simplicity of cellular algebras. {\it J.
Algebra} 267:323--341.\label{C1}

\item Cao, Y. Z., Zhu, P. (2006). Cyclotomic blob algebras and its representation theory. {\it J. Pure and applied Algebra} 204:666--695.\label{C03}

\item Cline, E., Parshall, B. J., Scott, L. L. (1988). Finite dimensional
algebras and highest weight categories. {\it J. Reine Angew. Math.}
391:85--99.\label{CPS88}

\item Fan, C. K., Green, R. M. (1999). On the affine Temperley-Lieb algebras. {\it
J. London Math. Soc.} 60:366--380.\label{FG99}

\item Freyd, P. J., Yetter, D. N. (1989). Braided compact closed categories with applications to low
dimensional topology. {\it Adv. Math.} 77:156--182.\label{FY89}

\item Graham, J. (1995). Modular Representations of Hecke Algebras and Related Algebras. Ph.D. thesis,
University of Sydney.\label{G95}

\item Graham, J., Lehrer, G. (1996). Cellular algebras. {\it Invent.
Math.} 123:1--34.\label{GL96}

\item Graham, J., Lehrer, G. (1998). The representation theory of affine Temperley-Lieb
algebras. {\it L'Enseignment Math.} 44:173--218.\label{GL98}

\item Green, R. M. (1998). Generalized Temperley-Lieb algebras and
decorated tangles. {\it J. Knot. Theory Ramifications}
7:155--171.\label{G98}

\item Jones, V. F. R. (1983). Index for subfactors. {\it Invent. Math.} 72:1--25.\label{J83}

\item Jones, V. F. R. (1987). Hecke algebra representations of braid groups and link polynomials. {\it Ann. of Math.} 126(2):335--388.\label{J87}

\item Kauffman, L. H. (1990). An invariant of regular isotopy. {\it Trans.
Amer. Math. Soc.} 318:417--471.\label{K90}

\item K\"onig, S., Wang, J. C. (2008). Cyclotomic extensions of
diagram algebras. {\it Comm. Algebra} 36(5):1739--1757.

\item K\"onig, S., Xi, C. C. (1998). On the structure of cellular
algebras. {\it in} ``Algebras and Modules II'' (I. Reiten, S.
Smal{\o} and {\O}. Solberg, Eds.). {\it Canadian Mathematical
Society Conference Proceedings} 24:365--386.\label{KX98}

\item K\"onig, S., Xi, C. C. (1999a). When is a cellular algebra
quasi-hereditary? {\it Math. Ann.} 315:281--293.\label{KX99M}

 \item K\"onig, S., Xi, C. C. (1999b). On the number of cells of a cellular
 algebra.
 {\it Comm. Algebra} 27:5463--5470.\label{KX99C}

 \item K\"onig S., Xi, C. C. (1999c). Cellular algebras and quasi-hereditary algebras: a
 comparison. {\it Electron. Res. Announc. Amer. Math. Soc.} 5:71--75.\label{KX99E}

\item K\"onig S., Xi, C. C. (2001). A characteristic free approach to Brauer algebras. {\it Trans.
Amer. Math. Soc.} 353:1489--1505.\label{KX01}

\item Martin, P. P., Saleur, H. (1994). The blob algebra and the periodic
 Temperley-Lieb algebra. {\it Lett. Math. Phys.} 30:189--206.\label{MS94}

\item Rui, H. B., Xi, C. C. (2004). The representation theory of
Cyclotomic Temperley-Lieb algebras. {\it Comment. Math. Helv.}
79(2):427--450.\label{RX}

\item Rui, H. B., Xi, C. C., Yu, W. H. (2005). On the semi-simplicity of cyclotomic Temperley-Lieb algebras. {\it Michigan Math. J.} 53(1):83--96.\label{RXYp}

\item Temperley, H. N. V., Lieb, E. H. (1971). Relations between
percolation and coloring problems and other graph theoretical
problems associated with regular planar lattices: some exact results
for the percolation problem. {\it Proc. Roy. Soc. London Ser. A}
322:251--280.\label{TL71}

\item Westbury, B. W. (1995). The representation theory of the
Temperley-Lieb algebra. {\it Math. Z.} 219:539--565.\label{W95}

\item Xi, C. C. (1999). Partition algebras are cellular. {\it Compos. Math.} 119:99--109.\label{X99}

\item Xi, C. C. (2000). On the quasi-heredity of Birman-Wenzl
algebras. {\it Adv. Math.} 154:280--298.\label{X00}

\item Xi, C. C. (2002). Standardly stratified algebras and
cellular algebras. {\it Math. Proc. Cambridge Philos. Soc.}
133:37--53.\label{X02}

\end{enumerate}

\end{document}